# Needlet algorithms for estimation in inverse problems

## Gérard Kerkyacharian

*Laboratoire de Probabilités et Modèles Aléatoires and Université Paris X-Nanterre (MODALX), 175 rue du Chevaleret, 75013 Paris, France*
*e-mail:* `kerk@math.jussieu.fr`

## Pencho Petrushev

*Department of Mathematics, University of South Carolina, Columbia, SC 29208*
*e-mail:* `pencho@math.sc.edu`

## Dominique Picard and Thomas Willer

*Laboratoire de Probabilités et Modèles Aléatoires and Université Paris VII-Denis Diderot, 175 rue du Chevaleret, 75013 Paris, France*
*e-mail:* `picard@math.jussieu.fr`; `willer@math.jussieu.fr`

**Abstract:** We provide a new algorithm for the treatment of inverse problems which combines the traditional SVD inversion with an appropriate thresholding technique in a well chosen new basis. Our goal is to devise an inversion procedure which has the advantages of localization and multiscale analysis of wavelet representations without losing the stability and computability of the SVD decompositions. To this end we utilize the construction of localized frames (termed "needlets") built upon the SVD bases.

We consider two different situations: the "wavelet" scenario, where the needlets are assumed to behave similarly to true wavelets, and the "Jacobi-type" scenario, where we assume that the properties of the frame truly depend on the SVD basis at hand (hence on the operator). To illustrate each situation, we apply the estimation algorithm respectively to the deconvolution problem and to the Wicksell problem. In the latter case, where the SVD basis is a Jacobi polynomial basis, we show that our scheme is capable of achieving rates of convergence which are optimal in the $L_2$ case, we obtain interesting rates of convergence for other $L_p$ norms which are new (to the best of our knowledge) in the literature, and we also give a simulation study showing that the NEED-D estimator outperforms other standard algorithms in almost all situations.

## 1. Introduction

We consider the problem of recovering a function $f$ from a blurred (by a linear operator) and noisy version of $f$: $Y_\varepsilon = Kf + \varepsilon \dot{W}$. It is important to note that, in general, for a problem like this there exists a basis which is fully adapted to the problem, and as a consequence, the inversion remains stable; this is the Singular Value Decomposition (SVD) basis. The SVD basis, however, might be difficult to determine and handle numerically. Also, it might not be appropriate for accurate description of the solution with a small number of parameters. Furthermore, in many practical situations, the signal exhibits inhomogeneous regularity, and its local features are particularly interesting to recover. In such cases, other bases or frames (in particular, localized wavelet type bases) might be much more appropriate to represent the object at hand.

Our goal is to devise an inversion procedure which has the advantages of localization and multiscale analysis of wavelet representations without losing the stability and computability of the SVD decompositions. To this end we utilize the construction (due to Petrushev and his co-authors) of localized frames (termed "needlets") built upon particular bases - here the SVD bases. This construction uses a Calderón type decomposition combined with an appropriate quadrature (cubature) formula. It has the big advantage of producing frames which are close to wavelet bases in terms of dyadic properties and localization, but because of their compatibility with the SVD bases provide stable and easily computable schemes.

NEED-D is an algorithm combining the traditional SVD inversion with an appropriate thresholding technique in a well chosen new basis. It enables one to approximate the targeted functions with excellent rates of convergence for any $\mathbb{L}_p$ loss function, and over a wide range of Besov spaces.

Our main idea is by combining the thresholding algorithm with SVD-based frames to create an effective and practically feasible algorithm for solving the inverse problem described above. The properties of the localized frame to be constructed depend on the underlying SVD basis. We will consider two different behaviors, the first corresponds to a "wavelet" behavior in the sense that the properties of the system are equivalent (as far as we are concerned) to the properties of a true wavelet basis. This case typically arises in the deconvolution setting. In the second case, the properties of the frame may differ from wavelet bases and truly depend on the SVD basis at hand (hence on the operator $K$). We will explore in detail a case typically arising when the SVD basis is a Jacobi polynomial basis. It is illustrated by the Wicksell problem. We show that our scheme is capable of achieving rates of convergence which are optimal in the $L_2$ case (to the best of our knowledge, for the Wicksell problem this is the only case studied up to now). For other $L_p$ norms we obtain interesting rates of convergence, which are new in the literature.

We also give a simulation study for the Wicksell problem which shows that the NEED-D algorithm applied in combination with SVD based frames is valuable since it outperforms other standard algorithms in almost all situations.

The paper is organized in the following way: the second section introduces the model, the classical SVD methods, and the two basic examples considered in this paper, i.e. the deconvolution and Wicksell problems. The third section introduces the needlet construction, gives some basic properties of needlets and introduces the NEED-D algorithm. The fourth section explores its properties in the wavelet scenario. The main motivation for the NEED-D algorithm is given there after. The fifth section is devoted to the results in a Jacobi scenario. The sixth section is devoted to simulation results. The proofs of the main results from sections 4–5 are given in sections 7–8, respectively. The last section is an appendix which contains the definition and basic properties of the Jacobi needlets and the associated Besov spaces.

## 2. Inverse Models

Suppose $\mathbb{H}$ and $\mathbb{K}$ are two Hilbert spaces and let $K : \mathbb{H} \mapsto \mathbb{K}$ be a linear operator. The standard linear ill-posed inverse problem consists in recovering a good approximation $f_\varepsilon$ of the solution $f$ of

$$g = Kf \tag{1}$$

when only a perturbation $Y_\varepsilon$ of $g$ is observed. In this paper, we will consider the case when this perturbation is an additive stochastic white noise. Namely, we observe $Y_\varepsilon$ defined by the following identity:

$$Y_\varepsilon = Kf + \varepsilon \dot{W}, \tag{2}$$

where $\varepsilon$ is the amplitude of the noise. It is supposed to be a small parameter which tends to 0. The error will be measured in terms of this small parameter. Here $\dot{W}$ is a $\mathbb{K}$-white noise, i.e. for any $g, h \in \mathbb{K}$, $\xi(g) := (\dot{W}, g)_\mathbb{K}$, $\xi(h) := (\dot{W}, h)_\mathbb{K}$ form random Gaussian vectors (centered) with marginal variance $\|g\|_\mathbb{K}^2$, $\|h\|_\mathbb{K}^2$, and covariance $(g, h)_\mathbb{K}$ (with the obvious extension when one considers $k$ functions instead of 2).

Equation (2) means that for any $g \in \mathbb{K}$, we observe $Y_\varepsilon(g) := (Y_\varepsilon, g)_\mathbb{K} = (Kf, g)_\mathbb{K} + \varepsilon \xi(g)$, where $\xi(g) \sim N(0, \|g\|^2)$, and $Y_\varepsilon(g)$, $Y_\varepsilon(h)$ are independent random variables for orthogonal functions $g$ and $h$.

### 2.1. The SVD Method

Under the assumption that $K$ is compact, there exist two orthonormal bases (SVD bases) $(e_k)$ of $\mathbb{H}$ and $(g_k)$ of $\mathbb{K}$, and a sequence $(b_k)$, $b_k \to 0$ as $k \to \infty$, such that

$$K^*Ke_k = b_k^2 e_k, \quad KK^*g_k = b_k^2 g_k, \quad Ke_k = b_k g_k; \quad K^*g_k = b_k e_k$$

with $K^*$ being the adjoint operator of $K$.
The Singular Value Decomposition (SVD) of $K$

$$Kf = \sum_k b_k \langle f, e_k \rangle g_k$$

gives rise to approximations of the type

$$f_\varepsilon = \sum_{k=0}^N b_k^{-1} \langle Y_\varepsilon, g_k \rangle e_k,$$

where $N = N(\varepsilon)$ has to be properly selected. This SVD method is very attractive theoretically and can be shown to be asymptotically optimal in many situations (see Mathe and Pereverzev [23] together with their non linear counterparts Cavalier and Tsybakov [6], Cavalier et al [5], Tsybakov [33], Goldenschluger and Pereverzev [16], Efromovich and Koltchinskii [12]). It also has the big advantage of performing a quick and stable inversion of the operator. However, it has serious limitations: First, the SVD bases might be difficult to determine and handle numerically. Secondly, while these bases are fully adapted to describe the operator $K$, they might not be appropriate for accurate description of the solution with a small number of coefficients. Also in many practical situations, the signal has inhomogeneous regularitiy, and its local features are particularly interesting to recover. In such cases, other bases (in particular, localized wavelet type bases) are much more suitable for representation of the object at hand.

In the last ten years, various nonlinear methods have been developed, especially in the direct case with the objective of automatically adapting to the unknown smoothness and local singular behavior of the solution. In the direct case, one of the most attractive methods is probably wavelet thresholding, since it allies numerical simplicity to asymptotic optimality on a large variety of functional classes such as Besov or Sobolev spaces.

To apply this approach to inverse problems, Donoho [10] introduced a wavelet-like decomposition, specifically adapted to the operator $K$ (Wavelet-Vaguelette-Decomposition) and utilized a thresholding algorithm to this decomposition. In Abramovich and Silverman [1], this method was compared with the similar vaguelette-wavelet decomposition. Other wavelet schemes should be mentioned here, such as the ones from Antoniadis and Bigot [3], Antoniadis & al [4], Dicken and Maass [9], and especially for the deconvolution problem, Penski & Vidakovic [29], Fan & Koo [13], Kalifa & Mallat [19], Neelamani & al [27]. Later on Cohen et al [7] introduced an algorithm combining a Galerkin inversion with a thresholding algorithm.

The approach developed here was greatly influenced by these works.

### 2.1.1. Deconvolution

The deconvolution problem is probably one of the most famous inverse problems, giving rise to a great deal of investigations, specially in signal processing, and has an extensive bibliography. In the deconvolution problem, we consider the following operator: Let in this case $\mathbb{H} = \mathbb{K}$ be the set of square integrable periodic functions, with the standard $\mathbb{L}_2[0,1]$ norm, and consider

$$f \in \mathbb{H} \mapsto Kf = \int_0^1 \gamma(u-t)f(t)dt \in \mathbb{H}, \qquad (3)$$

where $\gamma$ is a known function in $\mathbb{H}$. It is generally assumed to be a regular function. A standard example is the box-car function which plays an important role in extending this model to image processing and especially to analysis of sequences of images.

In this case simple calculations show that the SVD bases $e_k$ and $g_k$ both coincide with the Fourier basis. The singular values correspond to the Fourier coefficients of the function $\gamma$:
$$b_k = \hat{\gamma}_k. \qquad (4)$$

### 2.1.2. Wicksell's problem

Another typical example is the following classical Wicksell's problem [34]. Suppose a population of spheres is embedded in a medium. The spheres have radii that may be assumed to be drawn independently from a density $f$. A random plane slice is taken through the medium and those spheres that are intersected by the plane furnish circles which radii are the points of observation $Y_1, \ldots, Y_n$. The unfolding problem is then to determine the density of the sphere radii from the observed circle radii. This problem also arises in medicine, where the spheres might be tumors in an animal's liver (see Nyshka et al [28]), as well as in numerous other contexts (biological, engineering, etc.) see for instance Cruz-Orive [8].

The difficulty of estimating the target function is well illustrated by figure 1. The Wicksell operator has a smoothing effect, thus the local variations of the target function become almost invisible in the case of observations corrupted by noise. (Also compare the blurred and noised observations in figure 6 to the target functions of figure 4.)

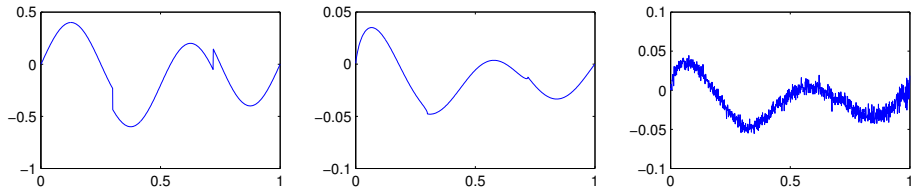

FIG 1. *Heavisine function, its image by the Wicksell operator without and with gaussian noise with $rsnr = 5$*

Following Wicksell [34] and Johnstone and Silverman [18], the Wicksell's problem corresponds to the following operator:
$$\mathbb{H} = \mathbb{L}_2([0,1], d\mu), \ d\mu(x) = (4x)^{-1} dx,$$
$$\mathbb{K} = \mathbb{L}_2([0,1], d\lambda), \ d\lambda(x) = 4\pi^{-1}(1-y^2)^{1/2} dy,$$

and
$$Kf(y) = \frac{\pi}{4} y(1-y^2)^{-1/2} \int_y^1 (x^2-y^2)^{-1/2} f(x) d\mu.$$

In this case, following [18], we have the following SVD bases:
$$e_k(x) = 4(k+1)^{1/2} x^2 P_k^{0,1}(2x^2 - 1)$$
$$g_k(y) = U_{2k+1}(y).$$

Here $P_k^{0,1}$ is the $k$th degree Jacobi polynomial of type $(0,1)$ and $U_k$ is the second type Chebishev polynomial of degree $k$. The singular values are

$$b_k = \frac{\pi}{16}(1+k)^{-1/2}. \tag{5}$$

In this article, in order to avoid some additional technicalities, we consider this problem in the white noise framework, which is simpler than the original problem described above in density terms.

## 3. General scheme for construction of frames (Needlets) and thresholding

Frames were introduced in the 1950's by Duffin and Schaeffer [11] as a means for studying nonharmonic Fourier series. These are redundant systems which behave like bases and allow for a lot of flexibility. Tight frame that are very close to orthonormal bases are particularly useful in signal and image processing.

In the following we present a general scheme for construction of frames due to Petrushev and his co-authors [26, 31, 30]. As will be shown this construction has the advantage of producing easily computable frame elements which are extremely well localized in all cases of interest. Following [26, 31, 30] we will call them "needlets".

Recall first the definition of a tight frame.

**Definition 1.** *Let $\mathbb{H}$ be a Hilbert space. A sequence $(\psi_n)$ in $\mathbb{H}$ is said to be a tight frame if*

$$\|f\|^2 = \sum_n |\langle f, \psi_n \rangle|^2 \quad \forall f \in \mathbb{H}.$$

Let $(\mathcal{Y}, \mu)$ be a measure space with $\mu$ a finite positive measure. Suppose we have the following decomposition

$$\mathbb{L}_2(\mathcal{Y}, \mu) = \bigoplus_{k=0}^{\infty} H_k,$$

where the $H_k$'s are finite dimensional spaces. For simplicity, we assume that $H_0$ is reduced to the constants.

Let $(e_i^k)_{i=1,\ldots,l_k}$ be an orthonormal basis of $H_k$. Then the orthogonal projector $L_k$ onto $H_k$ takes the form

$$L_k(f)(x) = \int_{\mathcal{Y}} f(y) L_k(x,y) d\mu(y), \quad \forall f \in \mathbb{L}_2(\mathcal{Y}, \mu),$$

where

$$L_k(x,y) = \sum_{i=1}^{l_k} e_i^k(x) \overline{e_i^k(y)}.$$

Note the obvious property of the orthogonal projectors:

$$\int_{\mathcal{Y}} L_k(x,y) L_m(y,z) d\mu(z) = \delta_{k,m} L_k(x,z). \tag{6}$$

The construction, inspired by the $\varphi$-transform of Frazier and Jawerth [15], consists of two main steps: (i) Calderón type decomposition and (ii) Discretization, which are described in the following two subsections.

### 3.1. Calderón type decomposition

Let $\varphi$ be a $C^\infty$ function supported in $[-1,1]$ such that $0 \leq \varphi(\xi) \leq 1$ and $\varphi(\xi) = 1$ if $|\xi| \leq \frac{1}{2}$. Define $a(\xi) \geq 0$ from

$$a^2(\xi) = \varphi(\xi/2) - \varphi(\xi) \geq 0.$$

Then

$$\sum_{j \geq 0} a^2(\xi/2^j) = 1, \quad \forall |\xi| \geq 1. \tag{7}$$

We now introduce the operator

$$\Lambda_j = \sum_{k \geq 0} a^2(k/2^j) L_k$$

and its associated kernel

$$\Lambda_j(x,y) = \sum_{k \geq 0} a^2(k/2^j) L_k(x,y) = \sum_{2^{j-1} < k < 2^{j+1}} a^2(k/2^j) L_k(x,y).$$

The operators $\Lambda_j$ provide a decomposition of $\mathbb{L}_2(\mathcal{Y}, \mu)$ which we record in the following proposition.

**Proposition 1.** *For all $f \in \mathbb{L}_2(\mathcal{Y}, \mu)$, we have*

$$f = L_0(f) + \sum_{j=0}^\infty \Lambda_j(f) \quad in \quad \mathbb{L}_2(\mathcal{Y}, \mu). \tag{8}$$

*Proof.* By the definition of $L_k$ and (7)

$$L_0 + \sum_{j=0}^J \Lambda_j = L_0 + \sum_{j=0}^J \sum_k a^2(k/2^j) L_k = \sum_k \varphi(k/2^{J+1}) L_k \tag{9}$$

and hence

$$\|f - L_0(f) - \sum_{j=0}^J \Lambda_j(f)\|^2 = \sum_{l \geq 2^{J+1}} \|L_l(f)\|^2 + \sum_{2^J \leq l < 2^{J+1}} \|L_l(f)(1 - \varphi(l/2^{J+1}))\|^2$$

$$\leq \sum_{l \geq 2^J} \|L_l(f)\|^2 \longrightarrow 0 \quad \text{as} \quad J \to \infty,$$

which completes the proof. □

### 3.2. Discretization

Let us define

$$\mathcal{K}_k = \bigoplus_{m=0}^k H_m.$$

We make two additional assumptions which will enable us to discretize decomposition (8) from Proposition 1:

(a)
$$f \in \mathcal{K}_k,\ g \in \mathcal{K}_l \Longrightarrow fg \in \mathcal{K}_{k+l}.$$

(b) *Quadrature formula*: For any $k \in \mathbb{N}$ there exists $\mathcal{X}_k$ a finite subset of $\mathcal{Y}$ and positive numbers $\lambda_\eta > 0$, $\eta \in \mathcal{X}_k$, such that

$$\int_\mathcal{Y} f d\mu = \sum_{\eta \in \mathcal{X}_k} \lambda_\eta f(\eta) \quad \forall f \in \mathcal{K}_k. \tag{10}$$

(Obviously, $\#\mathcal{X}_0 = 1$.)

We define
$$M_j(x,y) = \sum_k a(k/2^j) L_k(x,y) \quad \text{for} \quad j \geq 0. \tag{11}$$

Then as a consequence of (6), we have
$$\Lambda_j(x,y) = \int_\mathcal{Y} M_j(x,z) M_j(z,y) d\mu(z). \tag{12}$$

It is readily seen that $M_j(x,z) = \overline{M_j(z,x)}$ and

$$\forall\ x,\ z \mapsto M_j(x,z) \in \mathcal{K}_{2^{j+1}-1} \quad \text{and hence} \quad z \mapsto M_j(x,z) M_j(z,y) \in \mathcal{K}_{2^{j+2}-2}.$$

Now, by (10)
$$\Lambda_j(x,y) = \int_\mathcal{Y} M_j(x,z) M_j(z,y) d\mu(z) = \sum_{\eta \in \mathcal{X}_{2^{j+2}-2}} \lambda_\eta M_j(x,\eta) M_j(\eta,y),$$

which implies

$$\begin{aligned}
\Lambda_j f(x) &= \int_\mathcal{Y} \Lambda_j(x,y) f(y) d\mu(y) = \int_\mathcal{Y} \sum_{\eta \in \mathcal{X}_{2^{j+2}-2}} \lambda_\eta M_j(x,\eta) M_j(\eta,y) f(y) d\mu(y) \\
&= \sum_{\eta \in \mathcal{X}_{2^{j+2}-2}} \sqrt{\lambda_\eta} M_j(x,\eta) \int_\mathcal{Y} f(y) \overline{\sqrt{\lambda_\eta} M_j(y,\eta)} d\mu(y).
\end{aligned} \tag{13}$$

We are now prepared to introduce the desired frame. Let $\mathbb{Z}_j = \mathcal{X}_{2^{j+2}-2}$ for $j \geq 0$ and $\mathbb{Z}_{-1} = \mathcal{X}_0$. We define the frame elements (needlets) by

$$\psi_{j,\eta}(x) = \sqrt{\lambda_\eta} M_j(x,\eta), \quad \eta \in \mathbb{Z}_j,\ j \geq -1. \tag{14}$$

Notice that $\mathbb{Z}_{-1}$ consists of a single point and $\psi_0 = \psi_{-1,\eta}$, $\eta \in \mathbb{Z}_{-1}$, is the $\mathbb{L}_2$-normalized positive constant. Now (13) becomes

$$\Lambda_j f(x) = \sum_{\eta \in \mathbb{Z}_j} \langle f, \psi_{j,\eta} \rangle \psi_{j,\eta}(x). \tag{15}$$

**Proposition 2.** *The family $(\psi_{j,\eta})_{\eta \in \mathbb{Z}_j, j \geq -1}$ is a tight frame for $\mathbb{L}_2(\mathcal{Y},\mu)$.*

*Proof.* As
$$f = \lim_{J \longrightarrow \infty} L_0(f) + \sum_{j=0}^{J} \Lambda_j(f)$$

we have
$$\|f\|^2 = \lim_{J \longrightarrow \infty} \langle L_0(f), f \rangle + \sum_{j=0}^{J} \langle \Lambda_j(f), f \rangle.$$

But by (15)
$$\langle \Lambda_j(f), f \rangle = \sum_{\eta \in \mathbb{Z}_j} \langle f, \psi_{j,\eta} \rangle \langle \psi_{j,\eta}, f \rangle = \sum_{\eta \in \mathbb{Z}_j} |\langle f, \psi_{j,\eta} \rangle|^2, \quad j \geq 0,$$

and since $\psi_0$ is the normalized constant $\langle L_0(f), f \rangle = |\langle f, \psi_0 \rangle|^2$. Hence
$$\|f\|^2 = |\langle f, \psi_0 \rangle|^2 + \sum_{j \in \mathbb{N}_0, \ \eta \in \mathbb{Z}_j} |\langle f, \psi_{j,\eta} \rangle|^2,$$

which shows that $(\psi_{j,\eta})$ is a tight frame. □

### 3.3. Localization properties

The critical property of the frame construction above which makes it so attractive is the excellent localization of the frame elements (needlets) $(\psi_{j,\eta})$ in various settings of interest (see [25, 26, 31, 30]). The following figure (due to Paolo Baldi) is an illustration of this phenomenon. The rapidly oscillating function is the Legendre polynomial of degree $2^8$, whereas the localized one is a needlet constructed as explained above using Legendre polynomials of degree $\leq 2^8$ and centered approximately at zero. Its localization is remarkable taking into account that both functions are polynomials of the same degree.

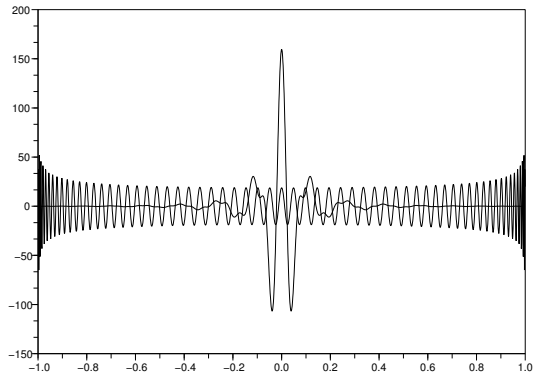

In the case of the unit sphere in $\mathbb{R}^{d+1}$, where $H_k$ are the spaces of spherical harmonics, the following localization property of the needlets is established in Narcowich, Petrushev and Ward [25, 26]: For any $k \in \mathbb{N}$ there exists a constant $C_k$ such that:
$$|\psi_{j\eta}(\xi)| \leq \frac{C_k 2^{dj/2}}{[1 + 2^j \arccos <\eta, \xi>]^k}.$$
In the case of Jacobi polynomials on $[-1, 1]$, the localization of the needlets proved in Petrushev, Xu [31] takes the form: For any $k \in \mathbb{N}$ there exists a constant $C_k$ such that
$$|\psi_{j\eta}(\cos\theta)| \leq \frac{C_k 2^{j/2}}{(1 + 2^j|\theta - \arccos\eta|)^k \sqrt{w_{\alpha,\beta}(2^j, \cos\theta)}}, \quad |\theta| \leq \pi,$$
where $w_{\alpha,\beta}(n, x) = (1 - x + n^{-2})^{\alpha+1/2}(1 + x + n^{-2})^{\beta+1/2}$ and $\alpha, \beta > -1/2$.

The almost exponential localization of the needlets and their semi-orthogonal structure allows to use them for characterization of spaces other than $L_2$, in particular the more general Triebel-Lizorkin and Besov spaces (see [26, 31]).

### *3.4. NEED-D algorithm: thresholding needlet coefficients*

We describe here the general idea of the method. The first step is to construct a needlet system (frame) $\{\psi_{j\eta} : \eta \in \mathbb{Z}_j, j \geq -1\}$ as described in section 3, where $H_k$ is simply the space spanned by the $k$-th vector $e_k$ of the SVD basis.

The needlet decomposition of any $f \in \mathbb{H}$ takes the form
$$f = \sum_j \sum_{\eta \in \mathbb{Z}_j} (f, \psi_{j\eta})_{\mathbb{H}} \psi_{j\eta}.$$

Using Parseval's identity, we have $\beta_{j\eta} = (f, \psi_{j\eta})_{\mathbb{H}} = \sum_i f_i \psi_{j\eta}^i$ with $f_i = (f, e_i)_{\mathbb{H}}$ and $\psi_{j\eta}^i = (\psi_{j\eta}, e_i)_{\mathbb{H}}$. If we put $Y_i = (Y_\varepsilon, g_i)_{\mathbb{K}}$, then
$$Y_i = (Kf, g_i)_{\mathbb{K}} + \varepsilon\xi_i = (f, K^*g_i)_{\mathbb{K}} + \varepsilon\xi_i = (\sum_j f_j e_j, K^*g_i)_{\mathbb{H}} + \varepsilon\xi_i = b_i f_i + \varepsilon\xi_i,$$
where $\xi_i = (\dot{W}, g_i)_{\mathbb{K}}$ form a sequence of centered Gaussian variables with variance 1. Thus
$$\hat{\beta}_{j\eta} = \sum_i \frac{Y_i}{b_i} \psi_{j\eta}^i$$
is an unbiased estimate of $\beta_{j\eta}$. Notice that from the needlet construction (see the previous section) it follows that the sum above is finite. More precisely, $\psi_{j\eta}^i \neq 0$ only for $2^{j-1} < i < 2^{j+1}$.

Let us consider the following estimate of $f$:
$$\hat{f} = \sum_{j=-1}^{J} \sum_{\eta \in \mathbb{Z}_j} t(\hat{\beta}_{j\eta}) \psi_{j\eta},$$

where $t$ is a thresholding operator defined by

$$t(\hat\beta_{j\eta}) = \hat\beta_{j\eta} I\{|\hat\beta_{j\eta}| \geq \kappa t_\varepsilon \sigma_j\} \quad \text{with} \tag{16}$$

$$t_\varepsilon = \varepsilon\sqrt{\log\frac{1}{\varepsilon}}. \tag{17}$$

Here $\kappa$ is a tuning parameter of the method which will be properly selected later on. Notice that the thresholding depends on the resolution level $j$ through the constant $\sigma_j$ which will also be specified later on, and the same with regard to the upper level of details $J$. Notice also that, in this paper, we concentrated on hard thresholding, whereas various other kind of thresholdings could be used, likely giving comparable results at least theoretically.

We will particularly focus on two situations (corresponding to the two examples discussed above). In the first case (see subsection 4), the needlets have very nice properties and behave exactly like wavelets. This is for instance the case of the deconvolution, where the SVD basis is the Fourier basis. However, more complicated problems e.g. the Wicksell's problem exhibit more delicate concentration properties for the needlets giving rise to different behaviors in terms of rates of convergence for the estimators.

## 4. NEED-D in wavelet scenario

In this section, we assume that the needlet system has the following properties: For any $1 \leq p < \infty$, there exist positive constants $c_p$, $C_p$, and $D_p$ such that

$$\text{Card } \mathbb{Z}_j \leq C 2^j, \tag{18}$$

$$c_p 2^{j(\frac{p}{2}-1)} \leq \|\psi_{j\eta}\|_p^p \leq C_p 2^{j(\frac{p}{2}-1)}, \tag{19}$$

$$\|\sum_{\eta\in\mathbb{Z}_j} u_\eta \psi_{j\eta}\|_p^p \leq D_p \sum_{\eta\in\mathbb{Z}_j} |u_\eta|^p \|\psi_{j\eta}\|_p^p, \text{ for any collection } (u_\eta). \tag{20}$$

We define the space $B^s_{\pi,r}$ as the collection of all functions $f$ with $f = \sum_{j\geq 0}\sum_{\eta\in\mathbb{Z}_j}\beta_{j\eta}\psi_{j\eta}$ such that

$$\|f\|_{B^s_{\pi,r}} := \|(2^{j[s+\frac{1}{2}-\frac{1}{\pi}]}\|(\beta_{j\eta})_{\eta\in\mathbb{Z}_j}\|_{l_\pi})_{j\geq 0}\|_{l_r} < \infty, \quad \text{and} \tag{21}$$

$$f \in B^s_{\pi,r}(M) \iff \|f\|_{B^s_{\pi,r}} \leq M. \tag{22}$$

**Theorem 1.** *Let $1 < p < \infty$, $2\nu + 1 > 0$, and*

$$\sigma_j^2 := \sup_\eta \sum_i [\frac{\psi^i_{j\eta}}{b_i}]^2 \leq C 2^{2j\nu}, \ \forall\, j \geq 0. \tag{23}$$

*Suppose $\kappa^2 \geq 16p$ and $2^J = [t_\varepsilon]^{\frac{-2}{2\nu+1}}$ with $t_\varepsilon$ as in (16).*

*Then for $f \in B^s_{\pi,r}(M)$ with $\pi \geq 1$, $s \geq 1/\pi$, $r \geq 1$ (with the restriction $r \leq \pi$ if $s = (\nu+\frac{1}{2})(\frac{p}{\pi}-1)$), we have*

$$\mathbb{E}\|\hat f - f\|_p^p \leq C \log(1/\varepsilon)^{p-1}[\varepsilon\sqrt{\log(1/\varepsilon)}]^{\mu p}, \tag{24}$$

*where*

$$\mu = \frac{s}{s+\nu+1/2}, \quad if \; s \geq (\nu + \frac{1}{2})(\frac{p}{\pi} - 1)$$

$$\mu = \frac{s - 1/\pi + 1/p}{s+\nu+1/2 - 1/\pi}, \quad if \; \frac{1}{\pi} \leq s < (\nu + \frac{1}{2})(\frac{p}{\pi} - 1).$$

The proof of this theorem is given in section 7.

*Remarks* :

1. These results are essentially minimax (see Willer [35]) up to logarithmic factors. We find back here the elbow, which was already observed in the direct problem, as well as in the deconvolution setting (see [17], for instance).
2. Condition (23) is essential in this problem. In the deconvolution case, the SVD basis is the *Fourier* basis and hence $\psi^i_{j\eta}$ are simply the Fourier coefficients of $\psi_{j\eta}$. Then assuming that we are in the so-called "regular" case ($b_k \sim k^{-\nu}$, for all $k$), it is easy to show that (23) is true for the needlet system as constructed above (see also the discussion in the following subsection). A similar remark can be made regarding conditions (19) and (20). In the deconvolution setting, the needlet construction is not strictly needed and, as is shown in Johnstone, Kerkyacharian, Picard, Raimondo[17], the periodized Meyer wavelet basis (see Meyer [24] and Mallat [22]) can replace the needlet construction. Condition (23) also holds in more general cases such as the box-car deconvolution, see [17], [20] where this algorithm is applied using Meyer's wavelets. ◇

### *4.1. Condition (23) and the needlet construction*

The following lines are intended to a posteriori motivate our decision to build upon the needlet construction. As was mentioned above condition (23) is very important for our algorithm. The proof will reveal that it is essential, since $\sigma^2_j$ is exactly the variance of our estimator of $\beta_{j\eta}$, so in a sense no other thresholding strategy can be better.

Let us now examine how condition (23) links the frame $(\psi_{j\eta})$ with the SVD basis $(e_k)$. To see this clearly let us suppose that $(\psi_{j\eta})$ is an arbitrary frame and let us place ourselves in the regular case:

$$b_i \sim i^{-\nu}$$

(this means that there exist two positive constants $c$ and $c'$ such that $c' i^{-\nu} \leq b_i \leq c i^{-\nu}$). If condition (23) holds true, we have

$$C 2^{2j\nu} \geq \sum_m \sum_{2^m \leq i \leq 2^{m+1}-1} [\frac{\psi^i_{j\eta}}{b_i}]^2.$$

Hence, $\forall\, m \geq j$,
$$\sum_{2^m \leq i \leq 2^{m+1}-1} [\psi_{j\eta}^i]^2 \leq c 2^{2\nu(j-m)}.$$

This means that the energy of $\psi_{j\eta}^i$ decays exponentially for $i \geq 2^j$, which reveals the role of the Littlewood Paley decomposition in the previous construction, replacing the exponential discrepancy by a cut-off.

The following proposition establishes a kind of converse property: The construction of needlet systems always implies that condition (23) is satisfied in the regular case.

**Proposition 3.** *If $(\psi_{j,\eta})$ is a frame such that $\{i : \psi_{j\eta}^i \neq 0\}$ is contained in a set $\{C_1 2^j, \ldots, C_2 2^j\}$, and $b_i \sim i^{-\nu}$, then*
$$\sigma_j^2 := \sum_i [\frac{\psi_{j\eta}^i}{b_i}]^2 \leq C 2^{2j\nu}.$$

*Proof.* Since the elements of an arbitrary frame are bounded in norm and $\psi_{j\eta}^i \neq 0$ only for $C_1 2^j \leq i \leq C_2 2^j$, we have
$$\sum_i [\frac{\psi_{j\eta}^i}{b_i}]^2 \leq C 2^{2j\nu} \|\psi_{j,\eta}\|^2 \leq C' 2^{2j\nu}.$$

□

## 5. NEED-D in a Jacobi-type case

Properties (19)–(20) are not necessarily valid for an arbitrary needlet system, since as mentioned above the localization properties of the frame elements depend on the initial underlying basis, and hence on the problem at hand. We will consider here a particular case motivated by Wicksell's problem.

We consider the space $\mathbb{H} = \mathbb{L}_2(I, d\gamma(x))$, where $I = [-1, 1]$, $d\gamma(x) = \omega_{\alpha,\beta}(x) dx$,
$$\omega_{\alpha,\beta}(x) = c_{\alpha,\beta}(1-x)^\alpha (1+x)^\beta, \quad \alpha, \beta > -1/2,$$

and $c_{\alpha,\beta}$ is selected so that $\int_I d\gamma_{\alpha,\beta}(x) = 1$. We will assume that $\alpha \geq \beta$ (otherwise we can interchange the roles of $\alpha$ and $\beta$).

Let $(P_k)_{k \geq 0}$ be the $\mathbb{L}_2(I, d\gamma(x))$ normalized Jacobi polynomials. We assume that the Jacobi polynomials appear as an SVD basis of the operator $K$. This is the case of Wicksell's problem, where $\beta = 0$, $\alpha = 1$, $b_k \sim k^{-1/2}$.

In the Jacobi case, the needlets have been introduced and studied in Petrushev and Xu [31]. See also the appendix, where the definition and some important properties of Jacobi needlets are given.

We will state our results in a more general setting, assuming that only a few conditions on the needlet system are valid. Note that these conditions are fulfilled by the needlet system (Jacobi needlets) constructed using the Jacobi polynomials $(P_k)_{k \geq 0}$. The proofs are given in the appendix.

We will consider two sets of conditions. The first one (which only depends on $\alpha$) is the following:

$$\text{Card } \mathbb{Z}_j \leq 2^j, \tag{25}$$

$$\sum_{\eta \in \mathbb{Z}_j} \|\psi_{j\eta}\|_p^p \leq C_p 2^{jp/2} \vee 2^{j(p-2)(1+\alpha)}, \quad \forall j, \; \forall p \neq 2 + \frac{1}{\alpha + 1/2}, \tag{26}$$

$$\|\sum_{\eta \in \mathcal{X}_j} \beta_\eta \psi_{j,\eta}\|_p \leq C (\sum_{\eta \in \mathcal{X}_j} |\beta_\eta|^p \|\psi_{j,\eta}\|_p^p)^{1/p}. \tag{27}$$

We define the space $\widetilde{B}^s_{\pi,r}$ as the collection of all functions $f$ on $[-1, 1]$ with representation

$$f = \sum_{j \geq -1} \sum_{\eta \in \mathbb{Z}_j} \beta_{j\eta} \psi_{j\eta}$$

such that

$$\|f\|_{\widetilde{B}^s_{\pi,r}} := \|(2^{js}(\sum |\beta_{j,\eta}|^\pi \|\psi_{j,\eta}\|_\pi^\pi)^{1/\pi})_{j \geq -1}\|_{l_r} < \infty, \quad \text{and} \tag{28}$$

$$f \in \widetilde{B}^s_{\pi,r}(M) \iff \|f\|_{\widetilde{B}^s_{\pi r}} \leq M. \tag{29}$$

In the Jacobi case, $\widetilde{B}^s_{\pi,r}$ is the Besov space defined in the Appendix as a space of approximation (not depending on the special needlet-frame).

**Theorem 2.** *Let $1 < p < \infty$ and $\alpha \geq \beta > -\frac{1}{2}$. Suppose*

$$t_\varepsilon = \varepsilon \sqrt{\log 1/\varepsilon} \quad \text{and} \quad 2^J = t_\varepsilon^{-\frac{2}{1+2\nu}}.$$

*Let $\kappa^2 \geq 16p[1 + 4\{(\frac{\alpha}{2} - \frac{\alpha+1}{p})_+ \vee (\frac{\beta}{2} - \frac{\beta+1}{p})_+\}]$ and assume that we are in the regular case, i.e.*

$$b_i \sim i^{-\nu}, \quad \nu > -\frac{1}{2}.$$

*Then for $f \in \widetilde{B}^s_{p,r}(M)$ with $s > [\frac{1}{2} - 2(\alpha + 1)(\frac{1}{2} - \frac{1}{p})]_+$, we have*

$$\mathbb{E}\|\hat{f} - f\|_p^p \leq C[\log(1/\varepsilon)]^{p-1}[\varepsilon\sqrt{\log(1/\varepsilon)}]^{\mu p},$$

*where*

(i) *if $p < 2 + \frac{1}{\alpha + 1/2}$, then*

$$\mu = \frac{s}{s + \nu + \frac{1}{2}};$$

(ii) *if $p > 2 + \frac{1}{\alpha + 1/2}$, then*

$$\mu = \frac{s}{s + \nu + (\alpha + 1)(1 - \frac{2}{p})}.$$

*Remarks* :

1. In the case $p < 2 + \frac{1}{\alpha+1/2}$, the rate obtained here is the usual one, and can be proved to be minimax (see [35]). The case $p > 2 + \frac{1}{\alpha+1/2}$ introduces a new rate of convergence.
2. Conditions (25)–(27) enabled us to estimate the rates of convergence of our scheme, whenever the index $\pi$ of the Besov space is the same as the index of the loss function ($p = \pi$). In the sequel, we will study the case where $p$ and $\pi$ are independently chosen. This requires, however, some additional assumptions. ◇

If in addition to properties (25)–(27), we now assume that the following conditions are fulfilled: For any $\eta \in \mathbb{Z}_j$, $j \geq 0$,

$$c2^{j(p-2)(\alpha+1)}k(\eta)^{-(p-2)(\alpha+1/2)} \leq \|\psi_{j\eta}\|_p^p \leq C2^{j(p-2)(\alpha+1)}k(\eta)^{-(p-2)(\alpha+1/2)},$$
$$k(\eta) < 2^{j-1}, \tag{30}$$
$$c2^{j(p-2)(\beta+1)}k'(\eta)^{-(p-2)(\beta+1/2)} \leq \|\psi_{j\eta}\|_p^p \leq C2^{j(p-2)(\beta+1)}k'(\eta)^{-(p-2)(\beta+1/2)},$$
$$k'(\eta) = 2^j - k(\eta) < 2^{j-1}, \tag{31}$$

where $k(\eta) \in \{1,\ldots,2^j\}$ is the index of $\eta \in \mathbb{Z}_j$. Here we assume that the points in $\mathbb{Z}_j$ are ordered so that $\eta_1 > \eta_2 > \cdots > \eta_{2^j}$. Note that in the case of Jacobi needlets $\mathbb{Z}_j$ consists of the zeros of the Jacobi polynomial $P_{2^j}^{\alpha,\beta}$ (see the appendix). In the following we will briefly write $k$ instead of $k(\eta)$ and $k'$ instead of $k'(\eta)$. Of course, (26) is now a consequence of conditions (30)–(31).

Observe the important fact that properties (30)–(31) are valid in the case of Jacobi Polynomials (see the appendix).

**Theorem 3.** *Let $1 < p < \infty$ and $\alpha \geq \beta > -\frac{1}{2}$. Suppose that conditions $(25) - (27)$ and $(30) - (31)$ are fulfilled. Let*

$$2^J = t_\varepsilon^{-\frac{2}{1+2\nu}} \quad and \quad \kappa^2 \geq 16p[1 + 4\{(\frac{\alpha}{2} - \frac{\alpha+1}{p})_+ \vee (\frac{\beta}{2} - \frac{\beta+1}{p})_+\}]$$

*and suppose that we are in the regular case, i.e.*

$$b_i \sim i^{-\nu}, \quad \nu > -\frac{1}{2}.$$

*Then for $f \in \widetilde{B}_{\pi,r}^s(M)$ with $s > \max_{\gamma \in \{\alpha,\beta\}}\{\frac{1}{2} - 2(\gamma+1)(\frac{1}{2} - \frac{1}{\pi}) \vee 2(\gamma+1)(\frac{1}{\pi} - \frac{1}{p}) \vee 0\}$, we have*

$$\mathbb{E}\|\hat{f} - f\|_p^p \leq C[\log(1/\varepsilon)]^{p-1+a}[\varepsilon\sqrt{\log(1/\varepsilon)}]^{\mu p}, \tag{32}$$

*where*

$$\mu = \min\{\mu(s), \mu(s,\alpha), \mu(s,\beta)\} \quad and \quad a = \max\{a(\alpha), a(\beta)\} \leq 2 \quad with$$

$$\mu(s) = \frac{s}{s+\nu+\frac{1}{2}},$$

$$\mu(s,\gamma) = \frac{s - 2(1+\gamma)(\frac{1}{\pi}-\frac{1}{p})}{s+\nu+2(1+\gamma)(\frac{1}{2}-\frac{1}{\pi})}$$

and,

$$a(\gamma) = \begin{cases} I\{\delta_p = 0\} & \text{if } [p-\pi][1-(p-2)(\gamma+1/2)] \geq 0, \\ \frac{(\gamma+\frac{1}{2})(\pi-p)}{(\pi-2)(\gamma+1/2)-1} + I\{\delta_s = 0\} & \text{if } [p-\pi][1-(p-2)(\gamma+1/2)] < 0, \end{cases}$$

with $\delta_p = 1 - (p-2)(\gamma+1/2)$ and $\delta_s = s[1-(p-2)(\gamma+1/2)] - p(2\nu+1)(\gamma+1)(\frac{1}{\pi}-\frac{1}{p})$.

The proofs of Theorems 2 and 3 are relegated to section 8.

*Remarks* :

1. Naturally, Theorem 2 follows by Theorem 3. We stated this two theorems separately because the hypotheses of Theorem 2 are less restrictive than the conditions in Theorem 3 and hence Theorem 2 potentially has wider range of application than Theorem 3.
2. It is interesting to notice that the convergence rates in (32) depend only on three distinctive regions for the parameters (which are actually present in Theorem 2, but hidden in the condition $\alpha \geq \beta$), which depends on a very subtle interrelation between the parameters $s$, $\alpha$, $\beta$, $p$, $\pi$.
3. It is also interesting to note that the usual rates of convergence obtained e.g. in the wavelet scenario are realized in the extreme case $\alpha = \beta = -\frac{1}{2}$.
$\diamond$

## 6. Simulation study

In this section we investigate the numerical performances of the NEED-D estimator in the context of the Wicksell problem described in section 2.1.2. We compare the results for simulated datasets to those obtained with several SVD methods.

### 6.1. The estimators

#### 6.1.1. Singular value decomposition estimators

With the notations introduced before, $f$ can be naturally estimated by the following linear estimator based on the singular value decomposition of operator $K$:

$$\hat{f} = \sum_i \lambda_i \frac{Y_i}{b_i} e_i,$$

where $(\lambda_i)_{i \in \mathbb{N}}$ is a deterministic filter.

In the simulations a first SVD estimator with projection weights was used:

$$\begin{cases} \lambda_i = 1 & \text{if } i \leq N, \\ \lambda_i = 0 & \text{if } i > N, \end{cases}$$

where the parameter $N$ was fitted for each setting so as to minimize the root mean square error ($RMSE$) of the estimator.

We also use the SVD estimator developed in Cavalier and Tsybakov [6], which is completely adaptive with a data driven choice of the filter and thus much more convenient than the former in practice. The values of $\lambda_i$ are constant in blocs $I_j = [\kappa_{j-1}, \kappa_j - 1]$ with limits $\kappa_0 = 1$ and $\kappa_J = N + 1$ determined further:

$$\begin{cases} \lambda_i = \left(1 - \frac{\sigma_j^2(1+\Delta_j^\gamma)}{\|\bar{Y}\|_{(j)}^2}\right)_+ & \text{if } i \in I_j, \; j = 1, \ldots J, \\ \lambda_i = 0 & \text{if } i > N, \end{cases}$$

where:

$$\bar{Y}_i = \frac{Y_i}{b_i}, \quad \|\bar{Y}\|_{(j)}^2 = \sum_{i \in I_j} \bar{Y}_i^2, \quad \sigma_j^2 = \varepsilon^2 \sum_{i \in I_j} b_i^{-2},$$

$$\Delta_j = \frac{\max_{i \in I_j} b_i^{-2}}{\sum_{i \in I_j} b_i^{-2}}, \quad 0 < \gamma < 1/2,$$

and we used the notation $x_+ = \max(0, x)$.

The blocks are determined by the following procedure. Let $\nu_\varepsilon \sim \max(5, \log\log(1/\varepsilon))$ and $\rho_\varepsilon = \frac{1}{\log(\nu_\varepsilon)}$, we define:

$$\begin{cases} \kappa_j = 1 & \text{if } j = 0, \\ \kappa_j = \nu_\varepsilon & \text{if } j = 1, \\ \kappa_j = \kappa_{j-1} + \lfloor \nu_\varepsilon \rho_\varepsilon (1+\rho_\varepsilon)^{j-1} \rfloor & \text{if } j = 2, \ldots, J, \end{cases}$$

where $J$ is large enough such that: $\kappa_J > \max\{m : \sum_{i=1}^m b_i^{-2} \leq \varepsilon^{-2} \rho_\varepsilon^{-3}\}$.

In the simulation settings considered further the value taken by $\kappa_J = N+1$ is too large compared to the level $n$ of the discretization resolution, thus the estimation was performed at the level $N_0 = \min\left(\frac{n}{2}, N\right)$ instead of $N$.

*6.1.2. Construction of the needlet basis*

Every needlet $\psi_{j,\eta_k}$ defined on $I = [-1, 1]$ is a linear combination of Jacobi polynomials as described in section 3, with weights depending on some filter $a$. This function is chosen as:

$$a(x) = \sqrt{\varphi(x/2) - \varphi(x)}, \quad \forall x \geq 0$$

where $\varphi(x) = I\{x < 0.5\} + P(x)I\{0.5 \leq x \leq 1\}$ and $P$ is a polynomial adjusted such that the corresponding needlet is sufficiently regular. In practice this choice seems to be slightly better than a $C^\infty$ filter with exponential shape.

The shape of $a$ is given by figure 2, and some examples of needlets are given in figure 3. Their amplitudes and supports fit automatically to the location of $\eta$: the needlets located near the edges of $I$ are much sharper than those located in the middle.

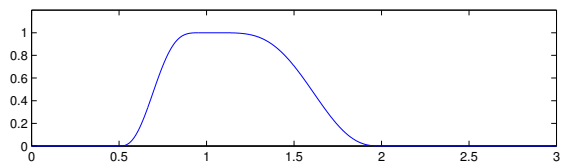

FIG 2. *Filter a with polynomial shape*

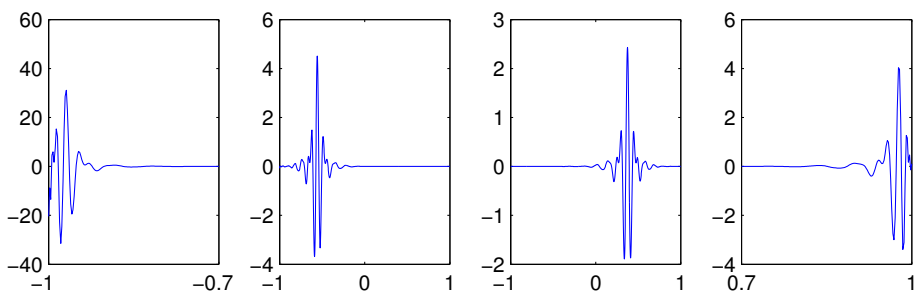

FIG 3. *Examples of needlets: $\psi_{7,\eta_{10}}$, $\psi_{7,\eta_{40}}$, $\psi_{7,\eta_{80}}$ and $\psi_{7,\eta_{120}}$ (from left to right)*

Finally NEED-D is performed by using the following basis $(\widetilde{\psi}_{j,\eta})$ adapted to the Wicksell problem:

$$\forall x \in [0,1], \quad \widetilde{\psi}_{j,\eta}(x) = 4\sqrt{2}x^2 \psi_{j,\eta}(2x^2 - 1).$$

With such a basis we have for all $i \in \mathbb{N}$:

$$\widetilde{\psi}^i_{j,\eta} = a(i/2^{j-1})P_i(\eta)\sqrt{b_{j,\eta}},$$

thus the estimated coefficients of $f$ in the frame are very easy to compute.

### *6.2. Parameters of the simulation*

We consider the four commonly used target functions $f$ represented in figure 4, and three levels of noise $\sigma$ corresponding to three values of the root signal to noise ratio of $K(f)$: $rsnr \in \{3, 5, 7\}$. The discretization resolution level is set to $n = 1024$, and the constant $\eta$ in the thresholds of NEED-D is set to $\eta = 0.75\sqrt{2}$.

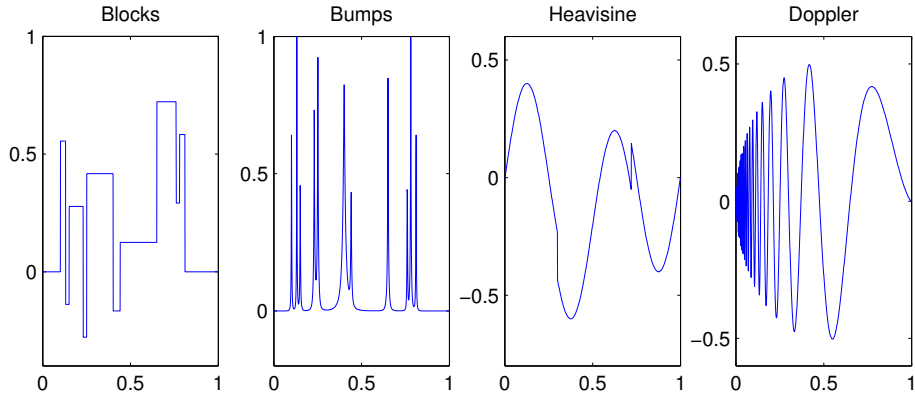

Fig 4. *Target functions*

The estimation error is evaluated by a Monte Carlo approximation of several $L_p(\mu)$ losses:

- $L1$ is computed as the average over 20 runs of $\frac{1}{n} \sum_{i=1}^{n} |f(\frac{i}{n}) - \hat{f}(\frac{i}{n})|/(\frac{4i}{n})$.
- $RMSE$ is computed as the average over 20 runs of $\sqrt{\frac{1}{n} \sum_{i=1}^{n} (f(\frac{i}{n}) - \hat{f}(\frac{i}{n}))^2/(\frac{4i}{n})}$.

In each run, the gaussian noise component is simulated independently of its values in the other runs.

### *6.3. Analysis of the results*

The performance of the non adaptive SVD estimator depends very strongly on the choice of $N$ (see figure 5). A large $N$ is needed in the case of small noise (first row of the figure) and in the case of very oscillating functions such as Doppler and Bumps. However even with this optimal *a posteriori* choice of $N$, the adaptive filter leads to better results than the non adaptive projection weights as shown in tables 1 and 2. Indeed the former is more adapted to the ill posed nature of the problem and to the variations of the noise, by adjusting over the singular values $(b_k)$ and the data $(y_k)$.

Moreover the NEED-D estimator generally outperforms both SVD estimators. As can be seen on figure 6, the differences are obvious in high noise for the Bumps and Doppler targets, where the SVD estimators are very noisy (in fact all the estimators happen to leave some noise unfiltered near the right edge

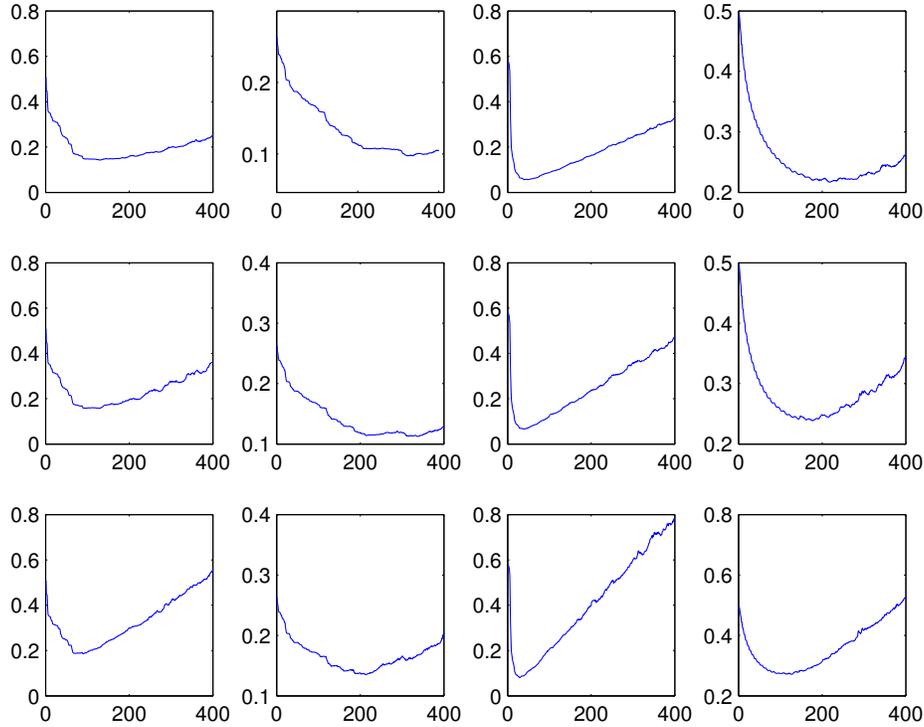

Fig 5. *Value of the mean square error of the non adaptive SVD estimator (y-axis) for each value of N (x-axis) for $rsnr = 7$ to $rsnr = 3$ (from top to bottom) and for the target function Blocks, Bumps, Heavisine and Doppler (from left to right)*

of the interval, which is given lesser importance by errors measured with the weight $\mu(x) = 1/(4x)$, for $x \in ]0,1]$.) This order of comparison is confirmed by the lower values of $L1$ and $RMSE$ for NEED-D than for SVD in all the settings (see tables 1 and 2).

|  | SVD | | | Adaptive SVD | | | NEED-D | | |
|---|---|---|---|---|---|---|---|---|---|
|  | low | med | high | low | med | high | low | med | high |
| **Blocks** | 0.0452 | 0.0495 | 0.0677 | 0.0399 | 0.0465 | 0.0591 | 0.0347 | 0.0404 | 0.0511 |
| **Bumps** | 0.0324 | 0.0388 | 0.0463 | 0.0258 | 0.0295 | 0.0361 | 0.0180 | 0.0206 | 0.0270 |
| **Heavisine** | 0.0257 | 0.0305 | 0.0402 | 0.0248 | 0.0299 | 0.0401 | 0.0205 | 0.0254 | 0.0321 |
| **Doppler** | 0.1032 | 0.1138 | 0.1307 | 0.1002 | 0.1085 | 0.1230 | 0.0858 | 0.0909 | 0.1007 |

TABLE 1
*Error L1 for each target, each noise level and each estimator*

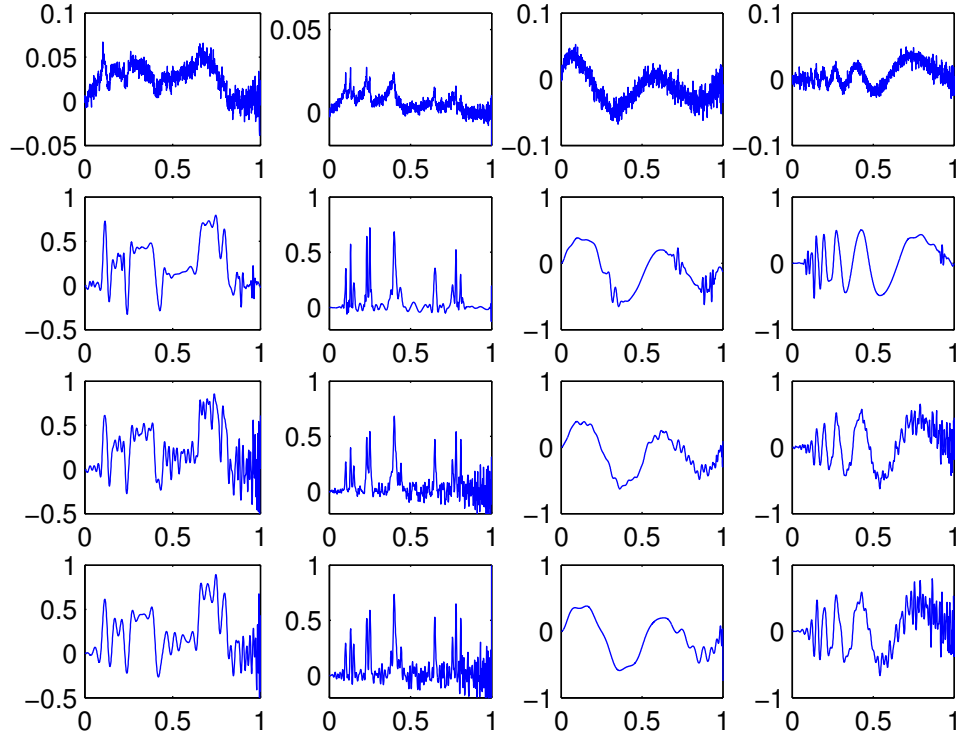

Fig 6. *From top to bottom: observed data, NEED-D estimator, adaptive SVD estimator and non adaptive SVD estimator for high noise (rsnr=3)*

## 7. Proof of Theorem 1

In this proof, $C$ will denote an absolute constant which may change from one line to the other.

|  | SVD | | | Adaptive SVD | | | NEED-D | | |
|---|---|---|---|---|---|---|---|---|---|
|  | low | med | high | low | med | high | low | med | high |
| **Blocks** | 0.0714 | 0.0790 | 0.0959 | 0.0665 | 0.0743 | 0.0900 | 0.0606 | 0.0673 | 0.0816 |
| **Bumps** | 0.0489 | 0.0577 | 0.0706 | 0.0453 | 0.0508 | 0.0617 | 0.0378 | 0.0416 | 0.0523 |
| **Heavisine** | 0.0278 | 0.0327 | 0.0422 | 0.0266 | 0.0317 | 0.0418 | 0.0235 | 0.0288 | 0.0379 |
| **Doppler** | 0.1092 | 0.1200 | 0.1378 | 0.1042 | 0.1114 | 0.1258 | 0.0969 | 0.0999 | 0.1071 |

TABLE 2
*Error L2 for each target, each noise level and each estimator*

First we have the following decomposition:

$$\mathbb{E}\|\hat{f} - f\|_p^p \le 2^{p-1}\{\mathbb{E}\|\sum_{j=-1}^{J}\sum_{\eta\in\mathbb{Z}_j}(t(\hat{\beta_{j\eta}}) - \beta_{j\eta})\psi_{j\eta}\|_p^p + \|\sum_{j>J}\sum_{\eta\in\mathbb{Z}_j}\beta_{j\eta}\psi_{j\eta}\|_p^p\}$$
$$=: I + II$$

The term $II$ is easy to analyse, as follows: Since $f$ belongs to $B_{\pi,r}^s(M)$, using standard embedding results (which in this case simply follows from direct comparisons between $l_q$ norms) we have that $f$ also belong to $B_{p,r}^{s-(\frac{1}{\pi}-\frac{1}{p})_+}(M')$, for some constant $M'$. Hence

$$\|\sum_{j>J}\sum_{\eta\in\mathbb{Z}_j}\beta_{j\eta}\psi_{j\eta}\|_p \le C2^{-J[s-(\frac{1}{\pi}-\frac{1}{p})_+]}.$$

Then we only need to verify that $\frac{s-(\frac{1}{\pi}-\frac{1}{p})_+}{\nu+1/2}$ is always larger that $\mu$, which is not difficult.

Bounding the term I is more involved. Using the triangular inequality together with Hölder inequality, and property (20) for the second line, we get

$$I \le 2^{p-1}J^{p-1}\sum_{j=-1}^{J}\mathbb{E}\|\sum_{\eta\in\mathbb{Z}_j}(t(\hat{\beta_{j\eta}}) - \beta_{j\eta})\psi_{j\eta}\|_p^p$$
$$\le 2^{p-1}J^{p-1}C\sum_{j=-1}^{J}\sum_{\eta\in\mathbb{Z}_j}\mathbb{E}|t(\hat{\beta_{j\eta}}) - \beta_{j\eta}|^p\|\psi_{j\eta}\|_p^p.$$

Now, we separate four cases:

$$\sum_{j=-1}^{J}\sum_{\eta\in\mathbb{Z}_j}\mathbb{E}|t(\hat{\beta_{j\eta}}) - \beta_{j\eta}|^p\|\psi_{j\eta}\|_p^p$$
$$= \sum_{j=-1}^{J}\sum_{\eta\in\mathbb{Z}_j}\mathbb{E}|t(\hat{\beta_{j\eta}}) - \beta_{j\eta}|^p\|\psi_{j\eta}\|_p^p\Big\{I\{|\hat{\beta_{j\eta}}| \ge \kappa t_\varepsilon\sigma_j\} + I\{|\hat{\beta_{j\eta}}| < \kappa t_\varepsilon\sigma_j\}\Big\}$$
$$\le \sum_{j=-1}^{J}\sum_{\eta\in\mathbb{Z}_j}\Big[\mathbb{E}|\hat{\beta_{j\eta}} - \beta_{j\eta}|^p\|\psi_{j\eta}\|_p^p I\{|\hat{\beta_{j\eta}}| \ge \kappa t_\varepsilon\sigma_j\}$$
$$\Big\{I\{|\beta_{j\eta}| \ge \frac{\kappa}{2}t_\varepsilon\sigma_j\} + I\{|\beta_{j\eta}| < \frac{\kappa}{2}t_\varepsilon\sigma_j\}\Big\}$$
$$+|\beta_{j\eta}|^p\|\psi_{j\eta}\|_p^p I\{|\hat{\beta_{j\eta}}| < \kappa t_\varepsilon\sigma_j\}\Big\{I\{|\beta_{j\eta}| \ge 2\kappa t_\varepsilon\sigma_j\} + I\{|\beta_{j\eta}| < 2\kappa t_\varepsilon\sigma_j\}\Big\}\Big]$$
$$\le: Bb + Bs + Sb + Ss.$$

If we notice that $\hat{\beta_{j\eta}} - \beta_{j\eta} = \sum_i \frac{Y_i - b_i f_i}{b_i}\psi_{j\eta}^i = \varepsilon\sum_i \xi_i\frac{\psi_{j\eta}^i}{b_i}$ is a gaussian random variable centered, and with variance $\varepsilon^2\sum_i[\frac{\psi_{j\eta}^i}{b_i}]^2$, we have using standard properties of the gaussian distribution, for any $q \geq 1$, if we recall that we set $\sigma_j^2 =: \sum_i [\frac{\psi_{jn}^i}{b_i}]^2 \leq C 2^{2j\nu}$, and denote by $s_q$ the $q$th absolute moment of the gaussian distribution when centered and with variance 1:

$$\mathbb{E}|\hat{\beta_{j\eta}} - \beta_{j\eta}|^q \leq s_q \sigma_j^q \varepsilon^q$$

$$\mathbb{P}\{|\hat{\beta_{j\eta}} - \beta_{j\eta}| \geq \frac{\kappa}{2} t_\varepsilon \sigma_j\} \leq 2\varepsilon^{\kappa^2/8}$$

Hence,

$$Bb \leq \sum_{j=-1}^{J} \sum_{\eta \in \mathbb{Z}_j} \sigma_j^p \varepsilon^p \|\psi_{j\eta}\|_p^p I\{|\beta_{j\eta}| \geq \frac{\kappa}{2} t_\varepsilon \sigma_j\}$$

$$Ss \leq \sum_{j=-1}^{J} \sum_{\eta \in \mathbb{Z}_j} |\beta_{j\eta}|^p \|\psi_{j\eta}\|_p^p I\{|\beta_{j\eta}| < 2\kappa t_\varepsilon \sigma_j\}.$$

And,

$$Bs \leq \sum_{j=-1}^{J} \sum_{\eta \in \mathbb{Z}_j} [\mathbb{E}|\hat{\beta_{j\eta}} - \beta_{j\eta}|^{2p}]^{1/2} [\mathbb{P}\{|\hat{\beta_{j\eta}} - \beta_{j\eta}|$$
$$\geq \frac{\kappa}{2} t_\varepsilon \sigma_j\}]^{1/2} \|\psi_{j\eta}\|_p^p I\{|\beta_{j\eta}| < \frac{\kappa}{2} t_\varepsilon \sigma_j\}$$
$$\leq \sum_{j=-1}^{J} \sum_{\eta \in \mathbb{Z}_j} s_{2p}^{1/2} \sigma_j^p \varepsilon^p 2^{1/2} \varepsilon^{\kappa^2/16} \|\psi_{j\eta}\|_p^p I\{|\beta_{j\eta}| < \frac{\kappa}{2} t_\varepsilon \sigma_j\}$$
$$\leq C \sum_{j=-1}^{J} 2^{jp(\nu+\frac{1}{2})} \varepsilon^p \varepsilon^{\kappa^2/16} \leq C \varepsilon^{\kappa^2/16}.$$

Now, if we remark that the $\beta_{j\eta}$ are necessarily all bounded by some constant (depending on $M$) since $f$ belongs to $B_{\pi,r}^s(M)$, and using (19),

$$Sb \leq \sum_{j=-1}^{J} \sum_{\eta \in \mathbb{Z}_j} |\beta_{j\eta}|^p \|\psi_{j\eta}\|_p^p \mathbb{P}\{|\hat{\beta_{j\eta}} - \beta_{j\eta}| \geq 2\kappa t_\varepsilon \sigma_j\} I\{|\beta_{j\eta}| \geq 2\kappa t_\varepsilon \sigma_j\}$$
$$\leq \sum_{j=-1}^{J} \sum_{\eta \in \mathbb{Z}_j} |\beta_{j\eta}|^p \|\psi_{j\eta}\|_p^p 2\varepsilon^{\kappa^2/8} I\{|\beta_{j\eta}| \geq 2\kappa t_\varepsilon t_\varepsilon \sigma_j\}$$
$$\leq C \sum_{j=-1}^{J} 2^{j\frac{p}{2}} \varepsilon^{\kappa^2/8} \leq C \varepsilon^{\frac{\kappa^2}{8} - \frac{p}{(2\nu+1)}}.$$

It is easy to check that in any cases if $\kappa^2 \geq 16p$ the terms $Bs$ and $Sb$ are smaller than the rates announced in the theorem.

If we recall that:
$$t_\varepsilon = \varepsilon \sqrt{\log \frac{1}{\varepsilon}}$$

We have using (19) and condition (23) for any $z \geq 0$:

$$Bb \leq C\varepsilon^p \sum_{j=-1}^{J} 2^{j(\nu p + \frac{p}{2} - 1)} \sum_{\eta \in \mathcal{Z}_j} I\{|\beta_{j\eta}| \geq \frac{\kappa}{2} t_\varepsilon \sigma_j\}$$

$$\leq C\varepsilon^p \sum_{j=-1}^{J} 2^{j(\nu p + \frac{p}{2} - 1)} \sum_{\eta \in \mathcal{Z}_j} |\beta_{j\eta}|^z [t_\varepsilon \sigma_j]^{-z}$$

$$\leq C t_\varepsilon^{p-z} \sum_{j=-1}^{J} 2^{j[\nu(p-z) + \frac{p}{2} - 1]} \sum_{\eta \in \mathcal{Z}_j} |\beta_{j\eta}|^z$$

Also, for any $p \geq z \geq 0$

$$Ss \leq C \sum_{j=-1}^{J} 2^{j(\frac{p}{2} - 1)} \sum_{\eta \in \mathcal{Z}_j} |\beta_{j\eta}|^z \sigma_j^{p-z} [t_\varepsilon]^{p-z}$$

$$\leq C[t_\varepsilon]^{p-z} \sum_{j=-1}^{J} 2^{j(\nu(p-z) + \frac{p}{2} - 1)} \sum_{\eta \in \mathcal{Z}_j} |\beta_{j\eta}|^z$$

So in both cases we have the same bound to investigate. We will write this bound on the following form (forgetting the constant):

$$I + II = t_\varepsilon^{p-z_1} [\sum_{j=-1}^{j_0} 2^{j[\nu(p-z_1) + \frac{p}{2} - 1]} \sum_{\eta \in \mathcal{Z}_j} |\beta_{j\eta}|^{z_1}]$$

$$+ t_\varepsilon^{p-z_2} [\sum_{j=j_0+1}^{J} 2^{j[\nu(p-z_2) + \frac{p}{2} - 1]} \sum_{\eta \in \mathcal{Z}_j} |\beta_{j\eta}|^{z_2}]$$

The constants $z_i$ and $j_0$ will be chosen depending on the cases, with the only constraint $p \geq z_i \geq 0$.

Notice first, that we only need to investigate the case $p \geq \pi$, since when $p \leq \pi$, $B_{\pi r}^s(M) \subset B_{pr}^s(M')$.

Let us first consider the case where $s \geq (\nu + \frac{1}{2})(\frac{p}{\pi} - 1)$, put

$$q = \frac{p(2\nu + 1)}{2(s + \nu) + 1}$$

and observe that on the considered domain, $q \leq \pi$ and $p > q$. In the sequel it will be useful to observe that we have $s = (\nu + \frac{1}{2})(\frac{p}{q} - 1)$. Now, taking $z_2 = \pi$, we get:

$$II \leq t_\varepsilon^{p-\pi} [\sum_{j=j_0+1}^{J} 2^{j[\nu(p-\pi) + \frac{p}{2} - 1]} \sum_{\eta \in \mathcal{Z}_j} |\beta_{j\eta}|^\pi]$$

Now, as

$$\frac{p}{2q} - \frac{1}{\pi} + \nu(\frac{p}{q} - 1) = s + \frac{1}{2} - \frac{1}{\pi}$$

and
$$\sum_{\eta \in \mathbb{Z}_j} |\beta_{j\eta}|^\pi = 2^{-j\pi(s+\frac{1}{2}-\frac{1}{\pi})}\tau_j^\pi$$

with $(\tau_j)_j \in l_r$ (this last thing is a consequence of the fact that $f \in B^s_{\pi,r}(M)$ and item (5)), we can write:

$$\begin{aligned} II &\leq t_\varepsilon^{p-\pi} \sum_{j=j_0+1} 2^{jp(1-\frac{\pi}{q})(\nu+\frac{1}{2})}\tau_j^\pi \\ &\leq C t_\varepsilon^{p-\pi} 2^{j_0 p(1-\frac{\pi}{q})(\nu+\frac{1}{2})} \end{aligned}$$

The last inequality is true for any $r \geq 1$ if $\pi > q$ and for $r \leq \pi$ if $\pi = q$. Notice that $\pi = q$ is equivalent to $s = (2\nu+1)(\frac{p}{2\pi} - \frac{1}{2})$. Now if we choose $j_0$ such that $2^{j_0 \frac{p}{q}(\nu+\frac{1}{2})} \sim t_\varepsilon^{-1}$ we get the bound

$$t_\varepsilon^{p-q}$$

which exactly gives the rate announced in the theorem for this case.

As for the first part of the sum (before $j_0$), we have, taking now $z_1 = \widetilde{q}$, with $\widetilde{q} \leq \pi$, so that $[\frac{1}{2^j}\sum_{\eta \in \mathbb{Z}_j}|\beta_{j\eta}|^{\widetilde{q}}]^{\frac{1}{\widetilde{q}}} \leq [\frac{1}{2^j}\sum_{\eta \in \mathbb{Z}_j}|\beta_{j\eta}|^\pi]^{\frac{1}{\pi}}$, and using again (7),

$$\begin{aligned} I &\leq t_\varepsilon^{p-\widetilde{q}}[\sum_{-1}^{j_0} 2^{j[\nu(p-\widetilde{q})+\frac{p}{2}-1]}\sum_{\eta \in \mathbb{Z}_j}|\beta_{j\eta}|^{\widetilde{q}}] \\ &\leq t_\varepsilon^{p-\widetilde{q}}[\sum_{-1}^{j_0} 2^{j[\nu(p-\widetilde{q})+\frac{p}{2}-\frac{\widetilde{q}}{\pi}]}[\sum_{\eta \in \mathbb{Z}_j}|\beta_{j\eta}|^\pi]^{\frac{\widetilde{q}}{\pi}}] \\ &\leq t_\varepsilon^{p-\widetilde{q}}\sum_{-1}^{j_0} 2^{j[(\nu+\frac{1}{2})p(1-\frac{\widetilde{q}}{q})]}\tau_j^{\widetilde{q}} \\ &\leq C t_\varepsilon^{p-\widetilde{q}} 2^{j_0[(\nu+\frac{1}{2})p(1-\frac{\widetilde{q}}{q})]} \\ &\leq C t_\varepsilon^{p-q} \end{aligned}$$

The last two lines are valid if $\widetilde{q}$ is chosen strictly smaller than $q$ (this is possible since $\pi \geq q$).

Let us now consider the case where $s < (2\nu+1)(\frac{p}{2\pi} - \frac{1}{2})$, and choose now

$$q = p\frac{\nu + \frac{1}{2} - \frac{1}{p}}{s + \nu + \frac{1}{2} - \frac{1}{\pi}}.$$

In such a way that we easily verify that $p - q = \frac{p(s+\frac{1}{p}-\frac{1}{\pi})}{s+\nu+\frac{1}{2}-\frac{1}{\pi}}$ and $q - \pi = \frac{(p-\pi)(\frac{1}{2}+\nu)-s\pi}{s+\nu+\frac{1}{2}-\frac{1}{\pi}} > 0$ Furthermore we also have $s + \frac{1}{2} - \frac{1}{\pi} = \frac{p}{2q} - \frac{1}{q} + \nu(\frac{p}{q} - 1)$.

Hence taking $z_1 = \pi$ and using again the fact that $f$ belongs to $B^s_{\pi,r}(M)$,

$$I \leq t_\varepsilon^{p-\pi}[\sum_{-1}^{j_0} 2^{j[\nu(p-\pi)+\frac{p}{2}-1]} \sum_{\eta \in \mathbb{Z}_j} |\beta_{j\eta}|^\pi]$$

$$\leq t_\varepsilon^{p-\pi} \sum_{-1}^{j_0} 2^{j[(\nu+\frac{1}{2}-\frac{1}{p})\frac{p}{q}(q-\pi)]} \tau_j^\pi$$

$$\leq C t_\varepsilon^{p-\pi} 2^{j_0[(\nu+\frac{1}{2}-\frac{1}{p})\frac{p}{q}(q-\pi)]}$$

This is true since $\nu + \frac{1}{2} - \frac{1}{p}$ is also strictly positive because of our constraints. If we now take $2^{j_0 \frac{p}{q}(\nu+\frac{1}{2}-\frac{1}{p})} \sim t_\varepsilon^{-1}$ we get the bound

$$t_\varepsilon^{p-q}$$

which is the rate announced in the theorem for this case.

Again, for $II$, we have, taking now $z_2 = \widetilde{q} > q(> \pi)$

$$II \leq t_\varepsilon^{p-\widetilde{q}} [\sum_{j=j_0+1}^{J} 2^{j[\nu(p-\widetilde{q})+\frac{p}{2}-1]} \sum_{\eta \in \mathbb{Z}_j} |\beta_{j\eta}|^{\widetilde{q}}]$$

$$\leq C t_\varepsilon^{p-\widetilde{q}} \sum_{j=j_0+1}^{J} 2^{j[(\nu+\frac{1}{2}-\frac{1}{p})\frac{p}{q}(q-\widetilde{q})]} z_j^{\frac{\widetilde{q}}{\pi}}$$

$$\leq C t_\varepsilon^{p-\widetilde{q}} 2^{j_0[(\nu+\frac{1}{2}-\frac{1}{p})\frac{p}{q}(q-\widetilde{q})]}$$

$$\leq C t_\varepsilon^{p-q}$$

## 8. Proof of the Theorems 2 and 3

The proof essentially follows the same steps as in the previous section. However, the following proposition will be helpful in the sequel.

**Proposition 4.** *Let us suppose that the following estimates are verified: Under the conditions (30) and (31), we have*

1.
$$\pi \geq p \Rightarrow (\sum_\eta |\beta_{j\eta}|^p \|\psi_{j,\eta}\|_p^p)^{1/p} \leq (\sum_\eta |\beta_{j\eta}|^\pi \|\psi_{j,\eta}\|_\pi^\pi)^{1/\pi}$$

2.
$$\pi < p \Rightarrow (\sum_{\eta, k(\eta) < 2^{j-1}} |\beta_{j\eta}|^p \|\psi_{j,\eta}\|_p^p)^{1/p}$$
$$\leq (\sum_{\eta, k(\eta) < 2^{j-1}} |\beta_{j\eta}|^\pi \|\psi_{j,\eta}\|_\pi^\pi)^{1/\pi} 2^{2j(\alpha+1)(1/\pi-1/p)}$$

$$\pi < p \Rightarrow (\sum_{\eta, k(\eta) \geq 2^{j-1}} |\beta_{j\eta}|^p \|\psi_{j,\eta}\|_p^p)^{1/p}$$
$$\leq (\sum_{\eta, k(\eta) \geq 2^{j-1}} |\beta_{j\eta}|^\pi \|\psi_{j,\eta}\|_\pi^\pi)^{1/\pi} 2^{2j(\beta+1)(1/\pi-1/p)}$$

*Proof.* • If $\pi \geq p$ clearly, because, Card$Z_j \leq 2^j$,

$$\left(\sum_{\eta \in Z_j} |\beta_{j\eta}|^p \|\psi_{j,\eta}\|_p^p\right)^{1/p} \leq 2^{j(1/p-1/\pi)} \left(\sum_{\eta \in Z_j} |\beta_{j\eta}|^\pi \|\psi_{j,\eta}\|_p^\pi\right)^{1/\pi}$$

But, using (30) and (31),

$$\pi \geq p \Rightarrow \|\psi_{j,\eta}\|_p \leq C \|\psi_{j,\eta}\|_\pi 2^{j(1/\pi-1/p)}.$$

So

$$\left(\sum_\eta |\beta_{j\eta}|^p \|\psi_{j,\eta}\|_p^p\right)^{1/p} \leq \left(\sum_\eta |\beta_{j\eta}|^\pi \|\psi_{j,\eta}\|_\pi^\pi\right)^{1/\pi}$$

• If $\pi \leq p$, clearly

$$\left(\sum_{\eta, k(\eta)<2^{j-1}} |\beta_{j\eta}|^p \|\psi_{j,\eta}\|_p^p\right)^{1/p} \leq \left(\sum_{\eta, k(\eta)<2^{j-1}} |\beta_{j\eta}|^\pi \|\psi_{j,\eta}\|_p^\pi\right)^{1/\pi}$$

But

$$\|\psi_{j,\eta}\|_p \leq C \|\psi_{j,\eta}\|_\pi 2^{j2(\alpha+1)(1/\pi-1/p)}, \ \forall \eta \ k(\eta) < 2^{j-1}$$

Hence,

$$\left(\sum_{\eta, k(\eta)<2^{j-1}} |\beta_{j\eta}|^p \|\psi_{j,\eta}\|_p^p\right)^{1/p} \leq \left(\sum_{\eta, k(\eta)<2^{j-1}} |\beta_{j\eta}|^\pi \|\psi_{j,\eta}\|_\pi^\pi\right)^{1/\pi} 2^{2j(\alpha+1)(1/\pi-1/p)}$$

The proof of the inequality with $\beta$ instead of $\alpha$ obviously is identical. □

Going back to the main stream of the proof, we first decompose:

$$\mathbb{E}\|\hat{f} - f\|_p^p \leq 2^{p-1}\{\mathbb{E}\|\sum_{j=-1}^J \sum_{\eta \in \mathbb{Z}_j}(t(\hat{\beta}_{j\eta}) - \beta_{j\eta})\psi_{j\eta}\|_p^p + \|\sum_{j>J}\sum_{\eta \in \mathbb{Z}_j} \beta_{j\eta}\psi_{j\eta}\|_p^p\}$$
$$=: \ I + II$$

• For $II$, using (27),

$$\|\sum_{j>J}\sum_{\eta \in \mathbb{Z}_j}\beta_{j\eta}\psi_{j\eta}\|_p^p \leq (\sum_{j>J}\|\sum_{\eta \in \mathbb{Z}_j}\beta_{j\eta}\psi_{j\eta}\|_p)^p \leq C[\sum_{j>J}(\sum_{\eta \in \mathbb{Z}_j}\|\beta_{j\eta}\psi_{j\eta}\|_p^p)^{1/p}]^p$$

If $\pi \geq p$, if we put $\delta = \frac{2}{1+2\nu}$, using $f \in \widetilde{B}_{p,r}^s(M)$,

$$II \leq C 2^{-Jsp} = C t_\varepsilon^{\delta sp}$$

If $\pi < p$, we decompose $II$ in the following way

$$\begin{aligned} II &\leq C\{[\sum_{j>J}(\sum_{\eta, k(\eta)<2^{j-1}}|\beta_{j\eta}|^p\|\psi_{j\eta}\|_p^p)^{1/p}]^p \\ &\quad + [\sum_{j>J}(\sum_{\eta, k(\eta)\geq 2^{j-1}}|\beta_{j\eta}|^p\|\psi_{j\eta}\|_p^p)^{1/p}]^p\} \\ &:= II(\alpha) + II(\beta). \end{aligned}$$

Now, using (4), and $f \in \widetilde{B}_{p,r}^s(M)$,

$$II(\alpha) \leq C[\sum_{j>J} 2^{-js}2^{j2(\alpha+1)(1/\pi-1/p)}]^p$$

If $s > 2(\alpha+1)(1/\pi - 1/p)$

$$II(\alpha) \leq C2^{-J(s-2(\alpha+1)(1/\pi-1/p))p} = Ct_\varepsilon^{\delta(s-2(\alpha+1)(1/\pi-1/p))p}.$$

The term $II(\beta)$ can be treated in the same way.
• For $I$

Using the triangular inequality together with Hölder inequality, and (27) for the second line, we get

$$\begin{aligned} I &\leq 2^{p-1}J^{p-1}\sum_{j=-1}^{J}\mathbb{E}\|\sum_{\eta\in\mathbb{Z}_j}(t(\hat{\beta_{j\eta}})-\beta_{j\eta})\psi_{j\eta}\|_p^p \\ &\leq 2^{p-1}J^{p-1}C\sum_{j=-1}^{J}\sum_{\eta\in\mathbb{Z}_j}\mathbb{E}|t(\hat{\beta_{j\eta}})-\beta_{j\eta}|^p\|\psi_{j\eta}\|_p^p \\ &\leq 2^{p-1}J^{p-1}C[I(\alpha)+I(\beta)] \end{aligned}$$

In the last line we separated as previously, in the sum $\eta \in \mathbb{Z}_j$, the indices $k(\eta) < 2^{j-1}$ and $k(\eta) \geq 2^{j-1}$. We will only investigate in the sequel $I(\alpha)$, since the argument for $I(\beta)$ goes in the same way.

Now, we separate four cases:

$$\sum_{j=-1}^{J}\sum_{\eta, k(\eta)<2^{j-1}}\mathbb{E}|t(\hat{\beta_{j\eta}})-\beta_{j\eta}|^p\|\psi_{j\eta}\|_p^p$$

$$= \sum_{j=-1}^{J}\sum_{\eta, k(\eta)<2^{j-1}}\mathbb{E}|t(\hat{\beta_{j\eta}})-\beta_{j\eta}|^p\|\psi_{j\eta}\|_p^p\Big\{I\{|\hat{\beta_{j\eta}}|\geq \kappa t_\varepsilon\sigma_j\}+I\{|\hat{\beta_{j\eta}}|<\kappa t_\varepsilon\sigma_j\}\Big\}$$

$$\leq \sum_{j=-1}^{J}\sum_{\eta, k(\eta)<2^{j-1}}\Big[\mathbb{E}|\hat{\beta_{j\eta}}-\beta_{j\eta}|^p\|\psi_{j\eta}\|_p^p I\{|\hat{\beta_{j\eta}}|\geq \kappa t_\varepsilon\sigma_j\}$$

$$\Big\{I\{|\beta_{j\eta}|\geq \frac{\kappa}{2}t_\varepsilon\sigma_j\}+I\{|\beta_{j\eta}|<\frac{\kappa}{2}t_\varepsilon\sigma_j\}\Big\}$$

$$+|\beta_{j\eta}|^p\|\psi_{j\eta}\|_p^p I\{|\hat{\beta_{j\eta}}|\geq \kappa t_\varepsilon\sigma_j\}\Big\{I\{|\beta_{j\eta}|\geq 2\kappa t_\varepsilon\sigma_j\}+I\{|\beta_{j\eta}|<2\kappa t_\varepsilon\sigma_j\}\Big\}\Big]$$

$$\leq: Bb + Bs + Sb + Ss$$

If we notice, as before, that $\hat{\beta_{j\eta}} - \beta_{j\eta} = \sum_i \frac{Y_i-b_i f_i}{b_i}\psi_{j\eta}^i = \varepsilon\sum_i \xi_i \frac{\psi_{j\eta}^i}{b_i}$ is a gaussian random variable centered, and with variance $\varepsilon^2\sum_i(\frac{\psi_{j\eta}^i}{b_i})^2$, we have using standard properties of the gaussian distribution, for any $q > 0$:

$$\mathbb{E}|\hat{\beta_{j\eta}} - \beta_{j\eta}|^q \le s_q[\varepsilon^2 \sum_i (\frac{\psi_{j\eta}^i}{b_i})^2]^{q/2} \le s_q \sigma_j^q \varepsilon^q \le C 2^{j\nu q} \varepsilon^q$$

$$\mathbb{P}\{|\hat{\beta_{j\eta}} - \beta_{j\eta}| \ge \frac{\kappa}{2}\varepsilon\sqrt{\log\frac{1}{\varepsilon}}\sigma_j\} \le c\varepsilon^{\kappa^2/8}$$

Hence,

$$\begin{aligned} Bb &\le \sum_{j=-1}^J \sum_{\eta, k(\eta)<2^{j-1}} \sigma_j^p \varepsilon^p \|\psi_{j\eta}\|_p^p I\{|\beta_{j\eta}| \ge \frac{\kappa}{2}\varepsilon\sqrt{\log\frac{1}{\varepsilon}}\sigma_j\} \\ Ss &\le \sum_{j=-1}^J \sum_{\eta, k(\eta)<2^{j-1}} |\beta_{j\eta}|^p \|\psi_{j\eta}\|_p^p I\{|\beta_{j\eta}| < 2\kappa\varepsilon\sqrt{\log\frac{1}{\varepsilon}}\sigma_j\} \end{aligned}$$

And,

$$\begin{aligned} Bs &\le \sum_{j=-1}^J \sum_{\eta, k(\eta)<2^{j-1}} [\mathbb{E}|\hat{\beta_{j\eta}} - \beta_{j\eta}|^{2p}]^{1/2} [\mathbb{P}\{|\hat{\beta_{j\eta}} - \beta_{j\eta}| \ge \frac{\kappa}{2}\varepsilon\sqrt{\log\frac{1}{\varepsilon}}\sigma_j\}]^{1/2} \\ &\quad \|\psi_{j\eta}\|_p^p I\{|\beta_{j\eta}| < \frac{\kappa}{2}\varepsilon\sqrt{\log\frac{1}{\varepsilon}}\sigma_j\} \\ &\le \sum_{j=-1}^J \sum_{\eta, k(\eta)<2^{j-1}} \sigma_{2p}^{1/2} \sigma_j^p \varepsilon^p c^{1/2} \varepsilon^{\kappa^2/16} \|\psi_{j\eta}\|_p^p I\{|\beta_{j\eta}| < \frac{\kappa}{2}\varepsilon\sqrt{\log\frac{1}{\varepsilon}}\sigma_j\} \\ &\le c'\varepsilon^p \varepsilon^{\kappa^2/16} \sum_{j=-1}^J 2^{jp\nu} \sum_{\eta \in \mathbb{Z}_j} \|\psi_{j\eta}\|_p^p \\ &\le c'\varepsilon^p \varepsilon^{\kappa^2/16} 2^{J(\nu p + (p/2) \vee (p-2)(1+\alpha))} \end{aligned}$$

using (26). Now, if we remark that the $\beta_{j\eta}$ are necessarily all bounded by some constant $M$, since $f \in \widetilde{B}_{p,r}^s(M)$,

$$\begin{aligned} Sb &\le \sum_{j=-1}^J \sum_{\eta, k(\eta)<2^{j-1}} |\beta_{j\eta}|^p \|\psi_{j\eta}\|_p^p \mathbb{P}\{|\hat{\beta_{j\eta}} - \beta_{j\eta}| \\ &\ge 2\kappa\varepsilon\sqrt{\log\frac{1}{\varepsilon}}\sigma_j\} I\{|\beta_{j\eta}| \ge 2\kappa\varepsilon\sqrt{\log\frac{1}{\varepsilon}}\sigma_j\} \\ &\le \sum_{j=-1}^J \sum_{\eta, k(\eta)<2^{j-1}} |\beta_{j\eta}|^p \|\psi_{j\eta}\|_p^p c\varepsilon^{\kappa^2/8} I\{|\beta_{j\eta}| \ge 2\kappa\varepsilon\sqrt{\log\frac{1}{\varepsilon}}\sigma_j\} \\ &\le c\varepsilon^{\kappa^2/8} \sum_{j=-1}^J \sum_{\eta \in \mathbb{Z}_j} \|\psi_{j\eta}\|_p^p \\ &\le c''\varepsilon^{\frac{\kappa^2}{8}} 2^{J(p/2 \vee (p-2)(1+\alpha))} \end{aligned}$$

It is easy to check that in any cases for $\kappa^2$ large enough, the terms $Bs$ and $Sb$ are smaller than the rates announced in the two theorems.

Now we focus on the bounds of Bb and Ss. Let $q \in [0, p]$, we always have:

$$\varepsilon^p \sum_{j=-1}^{J} \sum_{\eta, k(\eta) < 2^{j-1}} \sigma_j^p \|\psi_{j\eta}\|_p^p I\{\frac{|\beta_{j\eta}|}{\sigma_j} \geq \frac{\kappa}{2} t_\varepsilon\}$$

$$\leq \varepsilon^p \sum_{j=-1}^{J} \frac{\sum_{\eta, k(\eta) < 2^{j-1}} \sigma_j^p \|\psi_{j\eta}\|_p^p |\beta_{j\eta}|^q}{(\kappa \sigma_j t_\varepsilon/2)^q}$$

$$\leq \varepsilon^p (\kappa t_\varepsilon/2)^{-q} \sum_{j=-1}^{J} \sum_{\eta, k(\eta) < 2^{j-1}} \sigma_j^{p-q} \|\psi_{j\eta}\|_p^p |\beta_{j\eta}|^q$$

And

$$\sum_{j=-1}^{J} \sum_{\eta, k(\eta) < 2^{j-1}} |\beta_{j\eta}|^p \|\psi_{j\eta}\|_p^p I\{|\beta_{j\eta}| < 2\kappa t_\varepsilon \sigma_j\}$$

$$\leq \sum_{j=-1}^{J} (2\kappa t_\varepsilon \sigma_j)^{p-q} \sum_{\eta, k(\eta) < 2^{j-1}} |\beta_{j\eta}|^q \|\psi_{j\eta}\|_p^p$$

$$\leq (2\kappa t_\varepsilon)^{p-q} \sum_{j=-1}^{J} \sum_{\eta, k(\eta) < 2^{j-1}} \sigma_j^{p-q} \|\psi_{j\eta}\|_p^p |\beta_{j\eta}|^q.$$

So like in the wavelet scenario we have the same bound to investigate:

$$Bb + Ss \leq \sum_{j=-1}^{J} \sum_{\eta, k(\eta) < 2^{j-1}} (t_\varepsilon \sigma_j)^{p-q} \|\psi_{j\eta}\|_p^p |\beta_{j\eta}|^q,$$

then we use (30) and we separate as before the bound obtained in two terms $A$ and $B$ with some parameters $j_0$, $z_1$ and $z_2$ determined later, depending on the cases:

$$A := \sum_{j=-1}^{j_0} (t_\varepsilon \sigma_j)^{p-z_1} 2^{j(p-2)(\alpha+1)} \sum_{\eta, k(\eta) < 2^{j-1}} |\beta_{j\eta}|^{z_1} k^{-(p-2)(\alpha+1/2)}$$

$$B := \sum_{j=j_0+1}^{J} (t_\varepsilon \sigma_j)^{p-z_2} 2^{j(p-2)(\alpha+1)} \sum_{\eta, k(\eta) < 2^{j-1}} |\beta_{j\eta}|^{z_2} k^{-(p-2)(\alpha+1/2)}.$$

Let us first suppose that $p \leq \pi$ and $(p-2)(\alpha+1/2) \leq 1$, or that $p \geq \pi$ and $(p-2)(\alpha+1/2) \geq 1$.

Then we take $z_1 = 0$ and $z_2 = p$, and let us denote $\delta_p = 1 - (p-2)(\alpha + \frac{1}{2})$. We have:

$$\begin{aligned} A &= \sum_{j=-1}^{j_0} (t_\varepsilon \sigma_j)^p 2^{j(p-2)(\alpha+1)} \sum_{\eta, k(\eta) < 2^{j-1}} k^{-(p-2)(\alpha+1/2)} \\ &= \sum_{j=-1}^{j_0} (t_\varepsilon \sigma_j)^p 2^{j(p/2) \vee (p-2)(\alpha+1)} j^{I(\delta_p = 0)} \\ &\leq C(t_\varepsilon \sigma_{j_0})^p 2^{j_0(p/2) \vee (p-2)(\alpha+1)} (\log \frac{1}{\varepsilon})^{I(\delta_p = 0)}. \end{aligned}$$

And by treating $B$ as was done previously with the term $II(\alpha)$, we obtain:

$$\begin{aligned} B &= \sum_{j=j_0+1}^{J} 2^{j(p-2)(\alpha+1)} \sum_{\eta, k(\eta) < 2^{j-1}} |\beta_{j\eta}|^p k^{-(p-2)(\alpha+1/2)} \\ &\leq C 2^{-j_0 p[s - 2(\alpha+1)(\frac{1}{\pi} - \frac{1}{p})_+]}. \end{aligned}$$

So if $p \leq \pi$ and $(p-2)(\alpha+1/2) \leq 1$ we set $2^{j_0} = t_\varepsilon^{-1/[s+\nu+\frac{1}{2}]}$, which yields:

$$A + B \leq C t_\varepsilon^{p \frac{s}{s+\nu+\frac{1}{2}}} (\log \frac{1}{\varepsilon})^{I(\delta_p = 0)},$$

and if $p \geq \pi$ and $(p-2)(\alpha+1/2) \geq 1$ we take $2^{j_0} = t_\varepsilon^{-1/[s+\nu+(\alpha+1)(1-\frac{2}{\pi})]}$, which yields:

$$A + B \leq C t_\varepsilon^{p \frac{s - 2(\alpha+1)(\frac{1}{\pi} - \frac{1}{p})}{s+\nu+(\alpha+1)(1-\frac{2}{\pi})}} (\log \frac{1}{\varepsilon})^{I(\delta_p = 0)}.$$

*In the other cases: $p < \pi$ and $(p-2)(\alpha+1/2) > 1$, or $p > \pi$ and $(p-2)(\alpha+1/2) < 1$,* let us set $q = \frac{(p-2)(\alpha+1/2) - 1}{(\pi-2)(\alpha+1/2) - 1} \pi$, which satisfies:

$$p - q = \frac{2(\alpha+1)(\pi - p)}{(\pi-2)(\alpha+1/2) - 1}, \quad \text{and} \quad \pi - q = \frac{\pi(\alpha+1/2)(\pi - p)}{(\pi-2)(\alpha+1/2) - 1},$$

so $q \in ]0, p \wedge \pi[$ under the assumptions made above.

Let us bound the quantity $\sum_{\eta, k(\eta) < 2^{j-1}} |\beta_{j\eta}|^q k^{-(p-2)(\alpha+1/2)}$. We define:

$$\delta_1 = -\frac{q}{\pi}(\pi - 2)(\alpha + 1/2), \quad \text{and} \quad \delta_2 = -(p-2)(\alpha + 1/2) - \delta_1.$$

Using Hölder inequality, (30), and the fact that $f \in \widetilde{B}^s_{p,r}(M)$, we have:

$$\sum_{\eta, k(\eta) < 2^{j-1}} |\beta_{j\eta}|^q k^{-(p-2)(\alpha+1/2)}$$

$$= \sum_{\eta \in \mathbb{Z}_j} |\beta_{j\eta}|^q k^{\delta_1} k^{\delta_2}$$

$$\leq \Big[\sum_{\eta, k(\eta) < 2^{j-1}} |\beta_{j\eta}|^\pi k^{-(\pi-2)(\alpha+1/2)}\Big]^{\frac{q}{\pi}} \Big[\sum_{\eta, k(\eta) < 2^{j-1}} k^{\frac{\delta_2}{1-\frac{q}{\pi}}}\Big]^{1-\frac{q}{\pi}}$$

$$\leq C 2^{-jsq - j\frac{q}{\pi}(\pi-2)(\alpha+1)} \Big[\sum_{\eta, k(\eta) < 2^{j-1}} k^{\frac{\pi\delta_2}{\pi-q}}\Big]^{1-\frac{q}{\pi}}$$

$$= C 2^{-j(p-2)(\alpha+1)} 2^{j(-sq + \frac{p-q}{2})} j^{1-\frac{q}{\pi}}.$$

In the last line we used the fact that:

$$(p-2)(\alpha+1) - sq - \frac{q}{\pi}(\pi-2)(\alpha+1) = -sq + \frac{p-q}{2}, \quad \text{and} \quad \frac{\pi\delta_2}{\pi-q} = -1.$$

1. Let us assume that:

$$-sq + (p-q)(\nu + \frac{1}{2}) < 0,$$

i.e. that:

$$\frac{-s\pi[(p-2)(\alpha+1/2) - 1] + (\alpha+1)(\pi-p)(2\nu+1)}{(\pi-2)(\alpha+1/2) - 1} < 0.$$

Then we take $z_1 = 0$ and $z_2 = q$:

$$\begin{aligned}
A &= \sum_{j=-1}^{j_0} (t_\varepsilon \sigma_j)^p 2^{j(p-2)(\alpha+1)} \sum_{\eta, k(\eta) < 2^{j-1}} k^{-(p-2)(\alpha+1/2)} \\
&\leq (t_\varepsilon \sigma_{j_0})^p 2^{j_0(p/2) \vee (p-2)(\alpha+1)},
\end{aligned}$$

$$\begin{aligned}
B &= \sum_{j=j_0+1}^J (t_\varepsilon \sigma_j)^{p-q} 2^{j(p-2)(\alpha+1)} \sum_{\eta, k(\eta) < 2^{j-1}} |\beta_{j\eta}|^q k^{-(p-2)(\alpha+1/2)} \\
&\leq C\Big[\sum_{j=j_0+1}^J (t_\varepsilon \sigma_j)^{p-q} 2^{j(-sq + \frac{p-q}{2})}\Big] J^{1-\frac{q}{\pi}} \\
&\leq C(t_\varepsilon \sigma_{j_0})^{p-q} 2^{j_0(-sq + \frac{p-q}{2})} (\log \frac{1}{\varepsilon})^{1-\frac{q}{\pi}}.
\end{aligned}$$

If $(p-2)(\alpha+1/2) > 1$ we take $2^{j_0} = t_\varepsilon^{-1/[s+\nu+(\alpha+1)(1-\frac{2}{\pi})]}$, which yields:

$$A + B \leq C t_\varepsilon^{p \frac{s - 2(\alpha+1)(\frac{1}{\pi} - \frac{1}{p})}{s+\nu+(\alpha+1)(1-\frac{2}{\pi})}} (\log \frac{1}{\varepsilon})^{1-\frac{q}{\pi}},$$

and if $(p-2)(\alpha+1/2) < 1$ we take $2^{j_0} = t_\varepsilon^{-1/[s+\nu+\frac{1}{2}]}$, which yields:

$$A + B \leq Ct_\varepsilon^{p\frac{s}{s+\nu+\frac{1}{2}}}(\log\frac{1}{\varepsilon})^{1-\frac{q}{\pi}}.$$

Notice that, because of our conditions on $s$, we always have $j_0 \leq J$.

2. Let us now assume that:
$$\frac{-s\pi[(p-2)(\alpha+1/2)-1] + (\alpha+1)(\pi-p)(2\nu+1)}{(\pi-2)(\alpha+1/2)-1} > 0.$$

Then we take $z_1 = q$ and $z_2 = p$:

$$\begin{aligned}
A &\leq C[\sum_{j=-1}^{j_0}(t_\varepsilon\sigma_j)^{p-q}2^{j(-sq+\frac{p-q}{2})}]J^{1-\frac{q}{\pi}} \\
&\leq C(t_\varepsilon\sigma_{j_0})^{p-q}2^{j_0(-sq+\frac{p-q}{2})}(\log\frac{1}{\varepsilon})^{1-\frac{q}{\pi}},
\end{aligned}$$

and as before with the bias term $II(\alpha)$:

$$\begin{aligned}
B &\leq \sum_{j=j_0+1}^{J}2^{j(p-2)(\alpha+1)}\sum_{\eta,k(\eta)<2^{j-1}}|\beta_{j\eta}|^p k^{-(p-2)(\alpha+1/2)} \\
&\leq C2^{-j_0 p[s-2(\alpha+1)(\frac{1}{\pi}-\frac{1}{p})+]}.
\end{aligned}$$

If $\pi > p$ we take $2^{j_0} = t_\varepsilon^{-1/[s+\nu+\frac{1}{2}]}$, which yields:

$$A + B \leq Ct_\varepsilon^{p\frac{s}{s+\nu+\frac{1}{2}}}(\log\frac{1}{\varepsilon})^{1-\frac{q}{\pi}},$$

and if $\pi < p$ we take $2^{j_0} = t_\varepsilon^{-1/[s+\nu+(\alpha+1)(1-\frac{2}{\pi})]}$, which yields:

$$A + B \leq Ct_\varepsilon^{p\frac{s-2(\alpha+1)(\frac{1}{\pi}-\frac{1}{p})}{s+\nu+(\alpha+1)(1-\frac{2}{\pi})}}(\log\frac{1}{\varepsilon})^{1-\frac{q}{\pi}}.$$

3. Let us finally assume that:
$$-s\pi[(p-2)(\alpha+1/2)-1] + (\alpha+1)(\pi-p)(2\nu+1) = 0.$$

We take $z_1 = q$ and $z_2 = p$ as previously:

$$\begin{aligned}
A + B &\leq C\sum_{j=-1}^{j_0}t_\varepsilon^{p-q}j^{1-\frac{q}{\pi}} + C2^{-j_0 p[s-2(\alpha+1)(\frac{1}{\pi}-\frac{1}{p})+]} \\
&\leq Ct_\varepsilon^{p-q}(\log\frac{1}{\varepsilon})^{2-\frac{q}{\pi}} + C2^{-j_0 p[s-2(\alpha+1)(\frac{1}{\pi}-\frac{1}{p})+]}.
\end{aligned}$$

We proceed exactly like in the previous case, and we obtain the same rate whether $\pi \geq p$ or $\pi < p$:

$$A + B \leq Ct_\varepsilon^{p\frac{s}{s+\nu+\frac{1}{2}}}(\log\frac{1}{\varepsilon})^{2-\frac{q}{\pi}}.$$

We can sum up all the results for $Bb$ and $Ss$ (and thus on $I(\alpha)$) very simply:

if $2(\alpha+1)(\frac{1}{\pi}-\frac{1}{p}) < s$ and $s[1-(p-2)(\alpha+1/2)] \leq p(2\nu+1)(\alpha+1)(\frac{1}{\pi}-\frac{1}{p})$ then:
$$Bb + Ss \leq Ct_\varepsilon^{p\frac{s+2(\alpha+1)(\frac{1}{p}-\frac{1}{\pi})}{s+\nu+(\alpha+1)(1-\frac{2}{\pi})}}(\log\frac{1}{\varepsilon})^a,$$

if $s[1-(p-2)(\alpha+1/2)] > p(2\nu+1)(\alpha+1)(\frac{1}{\pi}-\frac{1}{p})$ then:
$$Bb + Ss \leq Ct_\varepsilon^{p\frac{s}{s+\nu+\frac{1}{2}}}(\log\frac{1}{\varepsilon})^a,$$

where the power of the log factor depends on the parameters:
$$a = \begin{cases} I\{\delta_p = 0\} & \text{if } [p-\pi][1-(p-2)(\alpha+1/2)] \geq 0, \\ \frac{(\alpha+\frac{1}{2})(\pi-p)}{(\pi-2)(\alpha+1/2)-1} + I\{\delta_s = 0\} & \text{if } [p-\pi][1-(p-2)(\alpha+1/2)] < 0, \end{cases}$$

with $\delta_p = 1-(p-2)(\alpha+1/2)$ and $\delta_s = s[1-(p-2)(\alpha+1/2)] - p(2\nu+1)(\alpha+1)(\frac{1}{\pi}-\frac{1}{p})$.

Note that the first term in the second case is bounded by 1, so we have $a \leq 2$ whatever the case.

## 9. Appendix: Needlets induced by Jacobi polynomials

The main references for this appendix are the two papers [31] of Petrushev and Xu, and [21] of Kyriazis, Petrushev and Xu.

### 9.1. Jacobi needlets: Definition and basic properties

In this section we apply the general scheme from §3 for the construction of Jacobi needlets. We begin by introducing some necessary notation. We denote $I = [-1,1]$ and $d\gamma_{\alpha,\beta}(x) = c_{\alpha,\beta}\omega_{\alpha,\beta}(x)dx$, where
$$\omega_{\alpha,\beta}(x) = (1-x)^\alpha(1+x)^\beta; \quad \alpha, \beta > -1/2,$$

and $c_{\alpha,\beta}$ is defined by $\int_I d\gamma_{\alpha,\beta}(x) = 1$. Assume $P^{\alpha,\beta}$ are the classical Jacobi polynomials (cf. e.g. [32]). Let $\Pi_k^{\alpha,\beta}$ be the Jacobi polynomial of degree $k$, normalized in $L_2(d\gamma_{\alpha\beta})$, i.e.
$$\int_I \Pi_k^{\alpha,\beta}\Pi_n^{\alpha,\beta} d\gamma_{\alpha,\beta} = \delta_{m,n}.$$

Let $a(\xi)$ be as in §3.1 with the additional condition: $a(\xi) > c > 0$ for $3/4 \leq \xi \leq 7/4$. Note that $\operatorname{supp} a \subset [1/2, 2]$. We define as in §3.1
$$\Lambda_j(x,y) = \sum_k a(k/2^j)\Pi_k^{\alpha,\beta}(x)\Pi_k^{\alpha,\beta}(y).$$

Let $\eta_\nu = \cos\theta_{j,\nu}$, $\nu = 1, 2, \ldots, 2^j$, be the zeros of the Jacobi polynomial $P_{2^j}$ ordered so that $\eta_1 > \eta_2 > \cdots > \eta_{2^j}$ and hence $0 < \theta_{j,1} < \theta_{j,2} < \cdots < \theta_{j,2^j} < \pi$. It is well known that (cf. [32])

$$\theta_{j,\nu} \sim \frac{\nu\pi}{2^j}. \tag{33}$$

We set

$$\mathcal{X}_j = \{\eta_\nu : \nu = 1, 2, \ldots, 2^j\}.$$

Let $\Pi_n$ denote the space of all polynomials of degree less than $n$. As is well known [32] the zeros of the Jacobi polynomial $P_{2^j}$ serve as knots of the Gaussian quadrature which is exact for all polynomials from $\Pi_{2^{j+1}-1}$, that is,

$$\int_I P \, d\gamma_{\alpha,\beta} = \sum_{\eta_\nu \in \mathcal{X}_j} b_{j,\eta_\nu} P(\eta_\nu), \quad \forall P \in \Pi_{2^{j+1}-1},$$

where the coefficients $b_{j,\eta_\nu} > 0$ are the Christoffel numbers [32] and $b_{j,\eta_\nu} \sim 2^{-j}\omega_{\alpha,\beta}(2^j;\eta_\nu)$ with

$$\omega_{\alpha,\beta}(2^j;x) := (1 - x + 2^{-2j})^{\alpha+1/2}(1 + x + 2^{-2j})^{\beta+1/2}.$$

We now define the Jacobi needlets by

$$\psi_{j,\eta_\nu}(x) = \sqrt{b_{j,\eta_\nu}} \Lambda_{2^j}(x, \eta_\nu), \quad \nu = 1, 2, \ldots, 2^j; \ j \geq 0,$$

and we set $\psi_0(x) = \psi_{-1,\eta}(x) = 1$, $\eta \in \mathcal{X}_{-1}$ with $\mathcal{X}_{-1}$ consisting of only one point $\eta = 0$. From Proposition 2, $(\psi_{j,\eta_\nu})$ is a tight frame of $\mathbb{L}_2(d\gamma_{\alpha\beta})$, i.e.

$$\|f\|_2^2 = \sum_{j \geq -1} \sum_{\eta \in \mathcal{X}_j} |\langle f, \psi_{j,\eta}\rangle|^2, \quad \forall f \in \mathbb{L}_2(d\gamma_{\alpha\beta}).$$

Hence

$$\|\psi_{j,\eta_\nu}\|_2 \leq 1. \tag{34}$$

Notice that (34) cannot be an equality since otherwise the needlet system $(\psi_{j,\eta_\nu})$ would be an orthonormal basis and this is impossible since

$$\sum_\nu \sqrt{b_{j,\eta_\nu}} \psi_{j,\eta_\nu} = \sum_\nu b_{j,\eta_\nu} L_{2^j}(x, \eta_\nu) = \int_I L_{2^j}(x, y) d\gamma(x) = 0.$$

We now recall the two main results from [31] which will be essential steps in our development.

**Theorem 4.** *For any $l \geq 1$ there exists a constant $C_l > 0$ such that*

$$|\psi_{j,\eta_\nu}(\cos\theta)| \leq C_l \frac{1}{\sqrt{\omega_{\alpha,\beta}(2^j,\cos\theta)}} \frac{2^{j/2}}{(1 + 2^j|\theta - \frac{\pi\nu}{2^j}|)^l}, \quad 0 \leq \theta \leq \pi. \tag{35}$$

Obviously

$$\omega_{\alpha,\beta}(2^j;\cos\theta) = (2\sin^2(\theta/2) + 2^{-2j})^{\alpha+1/2}(2\cos^2(\theta/2) + 2^{-2j})^{\beta+1/2}. \quad (36)$$

Therefore, $0 \leq \theta \leq \pi/2 \Longrightarrow \omega_{\alpha,\beta}(2^j,\cos\theta) \sim ((2^j\theta+1)^{2\alpha+1}2^{-j(2\alpha+1)})$ and hence

$$|\psi_{j,\eta_\nu}(\cos\theta)| \leq C_l \frac{2^{j(1+\alpha)}}{(1+2^j|\theta - \frac{\nu\pi}{2^j}|)^l} \frac{1}{(2^j\theta+1)^{\alpha+1/2}}, \quad 0 \leq \theta \leq \pi/2. \quad (37)$$

Similarly, from (36)

$$\pi/2 \leq \theta \leq \pi \Longrightarrow \omega_{\alpha,\beta}(2^j,\cos\theta) \sim (2^j(\pi-\theta)+1)^{2\beta+1}2^{-j(2\beta+1)}$$

and hence

$$|\psi_{j,\eta_\nu}(\cos\theta)| \leq C_l \frac{2^{j(1+\beta)}}{(1+2^j|\theta - \frac{\nu\pi}{2^j}|)^l} \frac{1}{(2^j(\pi-\theta)+1)^{\beta+1/2}}, \quad \pi/2 \leq \theta \leq \pi. \quad (38)$$

**Theorem 5.** *Let $0 < p \leq \infty$. Then*

$$\|\psi_{j,\eta_\nu}\|_p = \Big(\int_I |\psi_{j,\eta_\nu}(x)|^p d\gamma_{\alpha,\beta}\Big)^{1/p} \leq C_p \Big(\frac{2^j}{\omega_{\alpha,\beta}(2^j;\eta_\nu)}\Big)^{1/2-1/p}.$$

Using (33) and (36), we infer $\omega_{\alpha,\beta}(j;\eta_\nu) \sim 2^{-j(2\alpha+1)}\nu^{2\alpha+1}$ if $1 \leq \nu \leq 2^{j-1}$ and $\omega_{\alpha,\beta}(j;\eta_\nu) \sim 2^{-j(2\beta+1)}(2^j - \nu + 1)^{2\beta+1}$ if $2^{j-1} < \nu \leq 2^j$. Consequently,

$$1 \leq \nu \leq 2^{j-1} \Longrightarrow \quad \|\psi_{j,\eta_\nu}\|_p \leq C_p \left(\frac{2^{j(\alpha+1)}}{\nu^{\alpha+1/2}}\right)^{1-2/p}, \quad (39)$$

$$2^{j-1} < \nu \leq 2^j \Longrightarrow \quad \|\psi_{j,\eta_\nu}\|_p \leq C_p \left(\frac{2^{j(\beta+1)}}{(2^j - \nu + 1)^{\beta+1/2}}\right)^{1-2/p}. \quad (40)$$

**Corollary 1.** *Let $1 \leq p \leq \infty$ and $1/p + 1/q = 1$. Then*

$$\forall (j,\eta_\nu), \quad \|\psi_{j,\eta_\nu}\|_p \|\psi_{j,\eta_\nu}\|_q \leq C_p C_q$$

### 9.2. Estimation of the $\mathbb{L}_p$ norms of the needlets

Here we establish estimates (30)–(31) for the norms of the Jacobi needlets. In fact we only need to prove the lower bounds because the upper bounds are given above, see Theorem 5 and (39)–(40). We record these bounds in the following theorem. We want to express our thanks to Yuan Xu for communicating to us another proof of this result.

**Theorem 6.** $\forall\, 0 < p \leq \infty,\ \forall j \in \mathbb{N},$

$$\forall\, \nu = 1,\ldots,2^{j-1},\quad c_p \left(\frac{2^{j(\alpha+1)}}{\nu^{\alpha+1/2}}\right)^{1-2/p} \leq \|\psi_{j,\eta_\nu}\|_p \leq C_p \left(\frac{2^{j(\alpha+1)}}{\nu^{\alpha+1/2}}\right)^{1-2/p}$$

$\forall\, 2^{j-1} < \nu \leq 2^j$,

$$c_p \left(\frac{2^{j(\beta+1)}}{(1+(2^j-\nu))^{\beta+1/2}}\right)^{1-2/p} \leq \|\psi_{j,\eta_\nu}\|_p \leq C_p \left(\frac{2^{j(\beta+1)}}{(1+(2^j-\nu))^{\beta+1/2}}\right)^{1-2/p}$$

or equivalently

$$\|\psi_{j,\eta_\nu}\|_p \sim \left(\frac{2^j}{\omega_{\alpha,\beta}(2^j;\eta_\nu)}\right)^{1/2-1/p}. \tag{41}$$

A critical role in the proof of this theorem will play the following proposition.

**Proposition 5.** *Let $c^\diamond$ be an arbitrary positive constant. Then there exists a constant $c > 0$ such that*

$$\sum_{k=N}^{2N-1} [P_k^{\alpha,\beta}(\cos\theta)]^2 \geq c\omega_{\alpha,\beta}(N;\cos\theta)^{-1} \quad \text{for} \quad c^\diamond N^{-1} \leq \theta \leq \pi - c^\diamond N^{-1}, \quad N \geq 2. \tag{42}$$

*Proof.* The proof will rely on the well known asymptotic representation of Jacobi polynomials (sf. [32, Theorem 8.21.12, p. 195]): For any constants $c > 0$ and $\varepsilon > 0$

$$\left(\sin\frac{\theta}{2}\right)^\alpha \left(\cos\frac{\theta}{2}\right)^\beta P_n^{\alpha,\beta}(\cos\theta)$$
$$= N^{-\alpha} \frac{\Gamma(n+\alpha+1)}{n!} \left(\frac{\theta}{\sin\theta}\right)^{1/2} J_\alpha(N\theta) + \theta^{1/2}\mathcal{O}(n^{-3/2}) \tag{43}$$

for $cn^{-1} \leq \theta \leq \pi - \varepsilon$, where $N = n + (\alpha+\beta+1)/2$ and $J_\alpha$ is the Bessel function. Further, using the well known asymptotic identity

$$J_\alpha(z) = \left(\frac{2}{\pi z}\right)^{1/2} \cos(z+\gamma) + \mathcal{O}(z^{-3/2}), \quad z \to \infty \quad (\gamma = -\alpha\pi/2 - \pi/4), \tag{44}$$

one obtains (sf. [32, Theorem 8.21.13, p. 195])

$$P_n^{\alpha,\beta}(\cos\theta) = (\pi n)^{-1/2} \left(\sin\frac{\theta}{2}\right)^{-\alpha-1/2} \left(\cos\frac{\theta}{2}\right)^{-\beta-1/2} \{\cos(N\theta+\gamma)+(n\theta)^{-1}\mathcal{O}(1)\} \tag{45}$$

for $cn^{-1} \leq \theta \leq \pi - cn^{-1}$.

As is well known the Jacobi polynomials $P_k^{\alpha,\beta}$ and $P_{k+1}^{\alpha,\beta}$ have no common zeros and hence it suffices to prove (42) only for sufficiently large $N$. Also, $P_k^{\alpha,\beta}(-x) = (-1)^k P_k^{\beta,\alpha}(x)$ and therefore it suffices to prove (42) only in the case $c^\diamond N^{-1} \leq \theta \leq \pi/2$.

Denote by $F_N(\theta)$ the left-hand side quantity in (42). Then by (45), applied with $c = 1/2$, it follows that

$$F_N(\theta) \geq N^{-1}\theta^{-2\alpha-1} \sum_{k=N}^{2N-1} \left(c_1 \cos^2(k\theta + h(\theta)) - c_2(k\theta)^{-2}\right)$$
$$\geq c' N^{-1}\theta^{-2\alpha-1} \sum_{k=N}^{2N-1} \cos^2(k\theta + h(\theta)) - c''\theta^{-2\alpha-1}(N\theta)^{-2},$$

for $N^{-1} \leq \theta \leq \pi/2$, where $h(\theta) = (\alpha + \beta + 1)\theta/2 - \pi\alpha/2 - \pi/4$. It is easy to verify that for $\pi N^{-1} \leq \theta \leq \pi/2$

$$\sum_{k=N}^{2N-1} \cos^2(k\theta + h) = \frac{N}{2} + \frac{\sin N\theta}{2\sin\theta}\cos((3N-1)\theta + 2h) \geq \frac{N}{2}\left(1 - \frac{\pi}{2N\theta}\right) \geq \frac{N}{4}.$$

Therefore,

$$\begin{aligned} F_N(\theta) &\geq \theta^{-2\alpha-1}(c'/4 - c''(N\theta)^{-2}) \geq (c'/8)\theta^{-2\alpha-1} \\ &\geq c\omega_{\alpha,\beta}(N;\theta) \quad \text{for} \quad c^*N^{-1} \leq \theta \leq \pi/2, \end{aligned} \quad (46)$$

where $c^* = \max\{\pi, (8c''/c')^{1/2}\} > 0$.

It remains to establish (42) for $c^\diamond N^{-1} \leq \theta \leq c^*N^{-1}$. Denote $\delta = (\alpha+\beta+1)/2$. We now apply (43) with $c = c^\diamond$ and $\varepsilon = \pi/2$ to obtain using that $\Gamma(n+\alpha+1)/n! \sim n^\alpha$, $\sin\theta \sim \theta$, and (44)

$$\begin{aligned} \left[P_k^{\alpha,\beta}(\cos\theta)\right]^2 &\geq \theta^{-2\alpha}\left(c_1[J_\alpha((k+\delta)\theta)]^2 - c_2 k^{-3/2}\theta^{1/2}|J_\alpha((k+\delta)\theta)|\right) \\ &\geq c_1\theta^{-2\alpha}[J_\alpha((k+\delta)\theta)]^2 - c\theta^{-2\alpha}k^{-2}. \end{aligned}$$

Choose $\lambda$ so that $\theta = \frac{\lambda}{N}$ and $c^\diamond \leq \lambda \leq c^*$. Summing up above we get

$$\begin{aligned} F_N(\theta) &\geq c_1\theta^{-2\alpha}\sum_{k=N}^{2N-1}[J_\alpha((k+\delta)\theta)]^2 - c\theta^{-2\alpha}N^{-1} \\ &= c_1\theta^{-2\alpha}N\sum_{k=N}^{2N-1}\frac{1}{N}\left[J_\alpha\left(\frac{(k+\delta)\lambda}{N}\right)\right]^2 - c\theta^{-2\alpha}N^{-1} \\ &= c_1\theta^{-2\alpha}N\sum_{j=0}^{N-1}\frac{1}{N}\left[J_\alpha\left(\frac{j\lambda}{N} + \lambda + \frac{\delta\lambda}{N}\right)\right]^2 - c\theta^{-2\alpha}N^{-1}. \end{aligned}$$

Obviously, the last sum above involves only values of the Bessel function $J_\alpha(\theta)$ for $c^\diamond \leq \theta \leq c^*(2+\delta)$ and hence uniformly in $\lambda \in [c^\diamond, c^*]$

$$\left|\sum_{j=0}^{N-1}\frac{1}{N}\left[J_\alpha\left(\frac{j\lambda}{N} + \lambda + \frac{\delta\lambda}{N}\right)\right]^2 - \sum_{j=0}^{N-1}\frac{1}{N}\left[J_\alpha\left(\frac{j\lambda}{N} + \lambda\right)\right]^2\right| \longrightarrow 0, \quad N \longrightarrow \infty.$$

The second sum above can be viewed as a Riemann sum of the integral $\int_0^1 J_\alpha^2(\lambda(\theta+1))d\theta$, which is a continuous function of $\lambda \in [c^\diamond, c^*]$ and hence $\min_{\lambda \in [c^\diamond, c^*]}\int_0^1 J_\alpha^2(\lambda(\theta+1))d\theta \geq \widetilde{c} > 0$. Consequently, for sufficiently large $N$

$$\begin{aligned} F_N(\theta) &\geq \theta^{-2\alpha}(\widetilde{c}c_1 N/2 - cN^{-1}) \geq c\theta^{-2\alpha}N \\ &\geq c\omega_{\alpha,\beta}(N;\theta) \text{for} \quad c^\diamond N^{-1} \leq \theta \leq c^*N^{-1}. \end{aligned}$$

From this and (46) it follows that (42) holds for sufficiently large $N$ and this completes the proof of Proposition 5. □

**Proof of Theorem 6.** We first note that (sf. [32]) $\Pi_k^{\alpha,\beta}(x) \sim k^{1/2} P_k^{\alpha,\beta}(x)$ and hence

$$\|\psi_{j,\eta_\nu}\|_2^2 = b_{j,\nu} \sum_{2^{j-2} < k < 2^j} a^2(k/2^j)(\Pi_k^{\alpha,\beta}(\cos\theta_{j,\nu}))^2$$

$$\geq c\omega_{\alpha,\beta}(2^j; \eta_\nu) \sum_{2^{j-2} < k < 2^j} a^2(k/2^j)(P_k^{\alpha,\beta}(\cos\theta_{j,\nu}))^2$$

$$\geq c\omega_{\alpha,\beta}(2^j; \eta_\nu) \sum_{\frac{3}{4}2^j \leq k \leq \frac{7}{4}2^j} (P_k^{\alpha,\beta}(\cos\theta_{j,\nu}))^2.$$

Observe also that there exists a constant $c^\diamond > 0$ such that $c^\diamond/2^j \leq \theta_{j,\nu} \leq \pi - c^\diamond/2^j$, $\nu = 1, 2, \ldots, 2^j$. We now employ Proposition 5 and (34) to conclude that

$$0 < c \leq \|\psi_{j,\eta_\nu}\|_2 \leq 1. \tag{47}$$

We need to establish only the lower bound in Theorem 6. Recall first the upper bound from Theorem 5

$$\|\psi_{j,\eta_\nu}\|_p \leq C_p \Big(\frac{2^j}{\omega_{\alpha,\beta}(2^j; \eta_\nu)}\Big)^{1/2 - 1/p}, \quad 0 < p \leq \infty. \tag{48}$$

Suppose $2 < p < \infty$ and let $1/p + 1/q = 1$. By (47) and Hölder's inequality we have

$$0 < c \leq \|\psi_{j,\eta_\nu}\|_2^2 \leq \|\psi_{j,\eta_\nu}\|_p \|\psi_{j,\eta_\nu}\|_q \leq c\|\psi_{j,\eta_\nu}\|_p \Big(\frac{2^j}{\omega_{\alpha,\beta}(2^j; \eta_\nu)}\Big)^{1/2 - 1/q}$$

which yields

$$\|\psi_{j,\eta_\nu}\|_p \geq c\Big(\frac{2^j}{\omega_{\alpha,\beta}(2^j; \eta_\nu)}\Big)^{1/2 - 1/p}. \tag{49}$$

The case $p = \infty$ is similar. In the case $0 < p < 2$, we have using (47)

$$0 < c \leq \|\psi_{j,\eta_\nu}\|_2^2 \leq \|\psi_{j,\eta_\nu}\|_p^p \|\psi_{j,\eta_\nu}\|_\infty^{2-p} \leq c\|\psi_{j,\eta_\nu}\|_p^p \Big(\frac{2^j}{\omega_{\alpha,\beta}(2^j; \eta_\nu)}\Big)^{1-p/2},$$

which implies (49). The lower bound estimates in Theorem 6 follow by (49). □

### 9.3. Bounding for the norm of a linear combination of needlets

Our goal is to prove estimate (27), which we record in the following theorem:

**Theorem 7.** *Let $0 < p < \infty$. There exists a constant $A_p > 0$ such that for any collection of numbers $\{\lambda_\nu : \nu = 1, 2, \ldots, 2^j\}$, $j \geq 0$,*

$$\Big\|\sum_{\nu=1}^{2^j} \lambda_\nu \psi_{j,\eta_\nu}\Big\|_{\mathbb{L}_p(\gamma_{\alpha,\beta})}^p \leq A_p \sum_{\nu=1}^{2^j} |\lambda_\nu|^p \|\psi_{j,\eta_\nu}\|_{\mathbb{L}_p(\gamma_{\alpha,\beta})}^p. \tag{50}$$

and if $p = \infty$,

$$\|\sum_\nu \lambda_\nu \psi_{j,\eta_\nu}\|_{\mathbb{L}_\infty(\gamma_{\alpha,\beta})} \leq A_\infty \sup_{\nu=1}^{2^j} |\lambda_\nu| \|\psi_{j,\eta_\nu}\|_{\mathbb{L}_\infty(\gamma_{\alpha,\beta})}. \tag{51}$$

*Proof.* Let us begin with the simplest case $p = \infty$

$$\|\sum_\nu \lambda_\nu \psi_{j,\eta_\nu}\|_\infty \leq \sup_\nu |\lambda_\nu| \|\psi_{j,\eta_\nu}\|_\infty \|\sum_\nu \frac{\psi_{j,\eta_\nu}}{\|\psi_{j,\eta_\nu}\|_\infty}\|_\infty$$

We can conclude using the following lemma:

**Lemma 1.** $\exists C < \infty$ such that

$$\forall j \in \mathbb{N}, \quad \|\sum_\nu \frac{\psi_{j,\eta_\nu}}{\|\psi_{j,\eta_\nu}\|_\infty}\|_\infty \leq C$$

*Proof.* Using Theorem 4, and (49), we have

$$\frac{\psi_{j,\eta_\nu}(\cos\theta)}{\|\psi_{j,\eta_\nu}\|_\infty} \leq C_l \frac{1}{\sqrt{\omega_{\alpha,\beta}(2^j, \cos\theta)}} \frac{2^{j/2}}{(1 + |2^j\theta - \pi\nu|)^l} \frac{\sqrt{\omega_{\alpha,\beta}(2^j, \eta_\nu)}}{2^{j/2}}$$

and hence

$$\|\sum_\nu \frac{\psi_{j,\eta_\nu}}{\|\psi_{j,\eta_\nu}\|_\infty}\|_\infty \leq \sup_\theta C \sum_\nu \frac{\sqrt{\omega_{\alpha,\beta}(2^j, \eta_\nu)}}{\sqrt{\omega_{\alpha,\beta}(2^j, \cos\theta)}} \frac{1}{(1 + |2^j\theta - \pi\nu|)^l}$$

It remains to prove that this last quantity is bounded.

But using (36), one easily shows that (see [21]):

$$\frac{\sqrt{\omega_{\alpha,\beta}(2^j, \eta_\nu)}}{\sqrt{\omega_{\alpha,\beta}(2^j, \cos\theta)}} = \frac{\sqrt{\omega_{\alpha,\beta}(2^j, \cos\theta_\nu)}}{\sqrt{\omega_{\alpha,\beta}(2^j, \cos\theta)}} \leq (1 + |2^j\theta - \nu\pi|)^{1/2 + \alpha\vee\beta}$$

consequently,

$$\sup_\theta \sum_\nu \frac{\sqrt{\omega_{\alpha,\beta}(2^j, \eta_\nu)}}{\sqrt{\omega_{\alpha,\beta}(2^j, \cos\theta)}} \frac{1}{(1 + |2^j\theta - \pi\nu|)^l}$$

$$\leq \sup_\theta \sum_\nu \frac{1}{(1 + |2^j\theta - \pi\nu|)^{l-1/2-\alpha\vee\beta}} < \infty \tag{52}$$

for sufficiently large $l$. □

Let now $0 < p < \infty$. Consider the maximal operator

$$(M_s f)(x) = \sup_{J \ni x} \left(\frac{1}{|J|} \int_J |f(u)|^s du\right)^{1/s}, \quad s > 0,$$

where the supremum is taken over all intervals $J \subset [-1, 1]$ which contain $x$ and $|J|$ denotes the length of $J$. As elsewhere, let $\alpha \wedge \beta > -1/2$. It is well known that the weight $\omega_{\alpha,\beta}(x) = (1-x)^\alpha(1+x)^\beta$ on $[-1,1]$ belongs to the Muckenhoupt class $A_p$ with $p > 1$ if $\alpha \vee \beta \leq p-1$. Then in the weighted case the Fefferman-Stein maximal inequality (see [14] and [2]) can be stated as follows: If $1 < p, r < \infty$ and $\omega_{\alpha,\beta} \in A_p$, then for any sequence of functions $(f_k)$ on $[-1,1]$

$$\left\|\left(\sum_k (M_1 f_k)^r\right)^{1/r}\right\|_{\mathbb{L}_p(\gamma_{\alpha,\beta})} \leq C_{p,r} \left\|\left(\sum_k |f_k|^r\right)^{1/r}\right\|_{\mathbb{L}_p(\gamma_{\alpha,\beta})}.$$

Using that $M_1|f|^s = (M_s f)^s$ one easily infers from above that the following maximal inequality holds: If $0 < p, r < \infty$ and $0 < s < \min\{p, r, \frac{p}{\alpha \vee \beta + 1}\}$, then for any sequence of functions $(f_k)$ on $[-1,1]$

$$\left\|\left(\sum_k (M_s f_k)^r\right)^{1/r}\right\|_{\mathbb{L}_p(\gamma_{\alpha,\beta})} \leq C \left\|\left(\sum_k |f_k|^r\right)^{1/r}\right\|_{\mathbb{L}_p(\gamma_{\alpha,\beta})}. \qquad (53)$$

As in §9.1, let $\eta_\nu = \cos\theta_{j,\nu}$, $\nu = 1, 2, \ldots, 2^j$, be the zeros of the Jacobi polynomial $P_{2^j}^{\alpha,\beta}$. Set $\eta_0 = 1$, $\eta_{2^j+1} = -1$ and $\theta_{j,0} = 0$, $\theta_{j,2^j+1} = \pi$, respectively. Denote $I_\nu = [\frac{\eta_\nu + \eta_{\nu+1}}{2}, \frac{\eta_\nu + \eta_{\nu-1}}{2}]$ and put

$$H_\nu = h_\nu 1_{I_\nu} \quad \text{with} \quad h_\nu = \left(\frac{2^j}{\omega_{\alpha,\beta}(2^j; \eta_\nu)}\right)^{1/2},$$

where $1_{I_\nu}$ is the indicator function of $I_\nu$.

We next show that for any $s > 0$

$$|\psi_{j,\eta_\nu}(x)| \leq c(M_s H_\nu)(x), \quad x \in [-1,1], \quad \forall \nu = 1, 2, \ldots, 2^j, \; j \geq 0. \qquad (54)$$

Obviously, $(M_s 1_{I_\nu})(x) = 1_{I_\nu}(x)$ for $x \in I_\nu$. Let $x \in [-1,1] \setminus I_\nu$ and set $\cos\theta = x$, $\theta \in [0, \pi]$. Then

$$[(M_s 1_{I_\nu})(x)]^s \sim \frac{|I_\nu|}{|x - \eta_\nu|} \sim \frac{\eta_{\nu-1} - \eta_{\nu+1}}{|x - \eta_\nu|}$$

$$\sim \frac{\sin\frac{1}{2}(\theta_{j,\nu+1} - \theta_{j,\nu-1})\sin\frac{1}{2}(\theta_{j,\nu+1} + \theta_{j,\nu-1})}{\sin\frac{1}{2}|\theta - \theta_{j,\nu}|\sin\frac{1}{2}(\theta + \theta_{j,\nu})}$$

$$\sim \frac{2^{-j}\theta_{j,\nu}}{|\theta - \theta_{j,\nu}|(\theta + \theta_{j,\nu})}.$$

Using that $\theta_{j,\nu} \geq c_* 2^{-j}$ for some constant $c_* > 0$, one easily verifies the inequality

$$\frac{\theta_{j,\nu}}{\theta + \theta_{j,\nu}} \geq \frac{1}{(2 + c_*^{-1})(1 + 2^j|\theta - \theta_{j,\nu}|)}.$$

From above it follows that

$$(M_s 1_{I_\nu})(\cos\theta) \geq \frac{c}{(1 + 2^j|\theta - \theta_{j,\nu}|)^{2/s}}, \quad \theta \in [0, \pi],$$

which along with (35) (applied with $l \geq 2/s$) yields (54).

Combining (54) and (53) we get

$$\|\sum_{\nu=1}^{2^j} \lambda_\nu \psi_{j,\eta_\nu}\|_{\mathbb{L}_p(\gamma_{\alpha,\beta})}^p \leq c \sum_{\nu=1}^{2^j} |\lambda_\nu|^p \|H_\nu\|_{\mathbb{L}_p(\gamma_{\alpha,\beta})}^p. \tag{55}$$

Straightforward calculation show that $\|1_{I_\nu}\|_{\mathbb{L}_p(\gamma_{\alpha,\beta})} \sim \left(2^{-j}\omega_{\alpha,\beta}(2^j;\eta_\nu)\right)^{1/p}$ and hence, using Theorem 6,

$$\|H_\nu\|_{\mathbb{L}_p(\gamma_{\alpha,\beta})} \sim \left(\frac{2^j}{\omega_{\alpha,\beta}(2^j;\eta_\nu)}\right)^{1/2-1/p} \sim \|\psi_{j,\eta_\nu}\|_{\mathbb{L}_p(\gamma_{\alpha,\beta})}.$$

This coupled with (55) implies (57). □

**Corollary 2.** *We have, for $1 \leq p \leq \infty$,*

$$\begin{aligned}
\|\sum_{\nu=1}^{2^j} \langle f, \psi_{j,\eta_\nu}\rangle \psi_{j,\eta_\nu}\|_{\mathbb{L}_p(\gamma_{\alpha,\beta})} &\leq A_p (\sum_{\nu=1}^{2^j} |\langle f, \psi_{j,\eta_\nu}\rangle|^p \|\psi_{j,\eta_\nu}\|_{\mathbb{L}_p(\gamma_{\alpha,\beta})}^p)^{1/p} \\
&\leq A'_p \|f\|_{\mathbb{L}_p(\gamma_{\alpha,\beta})}
\end{aligned} \tag{56}$$

*Proof.* We have only to prove the righthand side inequality.

Let us first consider the case $p = 1$.

$$\sum_{\nu=1}^{2^j} |\langle f, \psi_{j,\eta_\nu}\rangle| \|\psi_{j,\eta_\nu}\|_{\mathbb{L}_1(\gamma_{\alpha,\beta})} \leq \sum_{\nu=1}^{2^j} \int |f|(x)|\psi_{j,\eta_\nu}(x)|d\mu(x)\|\psi_{j,\eta_\nu}\|_{\mathbb{L}_1(\gamma_{\alpha,\beta})}$$

$$= \int |f|(x)\{\sum_{\nu=1}^{2^j} \frac{|\psi_{j,\eta_\nu}(x)|}{\|\psi_{j,\eta_\nu}\|_\infty}\|\psi_{j,\eta_\nu}\|_\infty\|\psi_{j,\eta_\nu}\|_1\}d\mu(x) \leq CC_1 C_\infty \|f\|_1$$

using lemma 1 and corollary 1.

For $p = \infty$, we have

$$\sup_{1 \leq \nu \leq 2^j} |\langle f, \psi_{j,\eta_\nu}\rangle| \|\psi_{j,\eta_\nu}\|_{\mathbb{L}_\infty(\gamma_{\alpha,\beta})}$$

$$\leq \sup_{1 \leq \nu \leq 2^j} \int |f|(x)|\psi_{j,\eta_\nu}(x)|d\mu(x)\|\psi_{j,\eta_\nu}\|_{\mathbb{L}_\infty(\gamma_{\alpha,\beta})}$$

$$\leq \|f\|_\infty \sup_{1 \leq \nu \leq 2^j} \|\psi_{j,\eta_\nu}\|_1 \|\psi_{j,\eta_\nu}\|_\infty \leq C_1 C_\infty \|f\|_\infty$$

and again we have used corollary 1.

Let now $1 < p < \infty$ Using Hölder's inequality $(1/p + 1/q = 1)$ we get

$$|\langle f, \psi_{j,\eta_\nu}\rangle|^p \leq \{\int |f||\psi_{j,\eta_\nu}|^{1/p}|\psi_{j,\eta_\nu}|^{1/q}d\mu\}^p \leq \int |f|^p |\psi_{j,\eta_\nu}|d\mu (\int |\psi_{j,\eta_\nu}|d\mu)^{p/q}$$

$$\leq \int |f|^p |\psi_{j,\eta_\nu}|d\mu \|\psi_{j,\eta_\nu}\|_1^{p/q}$$

Hence

$$\sum_{\nu=1}^{2^j} |\langle f, \psi_{j,\eta_\nu}\rangle|^p \|\psi_{j,\eta_\nu}\|_p^p \leq \sum_{\nu=1}^{2^j} \int |f|^p |\psi_{j,\eta_\nu}| d\mu \|\psi_{j,\eta_\nu}\|_1^{p/q} \|\psi_{j,\eta_\nu}\|_p^p$$

$$= \int |f|^p \{\sum_{\nu=1}^{2^j} |\psi_{j,\eta_\nu}(x)| \|\psi_{j,\eta_\nu}\|_1^{p/q} \|\psi_{j,\eta_\nu}\|_p^p\} d\mu(x)$$

$$= \int |f|^p \{\sum_{\nu=1}^{2^j} \frac{|\psi_{j,\eta_\nu}(x)|}{\|\psi_{j,\eta_\nu}\|_\infty} \|\psi_{j,\eta_\nu}\|_\infty \|\psi_{j,\eta_\nu}\|_1^{p/q} \|\psi_{j,\eta_\nu}\|_p^p\} d\mu(x) \leq C C_\infty C_1^{p/q} C_p^p \|f\|_p^p$$

as by (39) and (40) we have

$$\|\psi_{j,\eta_\nu}\|_\infty \|\psi_{j,\eta_\nu}\|_1^{p/q} \|\psi_{j,\eta_\nu}\|_p^p \leq C_\infty C_1^{p/q} C_p^p$$

we have concluded using lemma 1.

$\square$

**Corollary 3.** *let*

$$\Lambda_j(f) = \sum_{\nu=1}^{2^j} \langle f, \psi_{j,\eta_\nu}\rangle \psi_{j,\eta_\nu}$$

$$\begin{aligned}\|\Lambda_j(f)\|_{\mathbb{L}_p(\gamma_{\alpha,\beta})} &\leq C_p (\sum_{\nu=1}^{2^j} |\langle f, \psi_{j,\eta_\nu}\rangle|^p \|\psi_{j,\eta_\nu}\|_{\mathbb{L}_p(\gamma_{\alpha,\beta})}^p)^{1/p} \\ &\leq C_p' \|\Lambda_{j-1}f + \Lambda_j f + \Lambda_{j+1}f\|_{\mathbb{L}_p(\gamma_{\alpha,\beta})}\end{aligned} \quad (57)$$

*(by convention $\Lambda_{-1}(f) = 0$)*

This claim is a simple consequence of the previous corollary, as

$$\forall \eta_n u \in \mathcal{X}_j, \quad \langle f, \psi_{j,\eta_\nu}\rangle = \langle \Lambda_{j-1}f + \Lambda_j f + \Lambda_{j+1}f, \psi_{j,\eta_\nu}\rangle$$

### 9.4. Besov Space

For $f \in \mathbb{L}_p(\gamma_{\alpha,\beta})$ let us define:

$$\forall n \in \mathbb{N}, \quad E_n(f,p) = \inf_{P \in \Pi_n} \|f - P\|_p$$

Clearly

$$\|f\|_p \geq E_0(f,p) \geq E_1(f,p)....$$

**Definition 2.** *Besov space $B^s_{p,q}$.*

Let $1 \leq p \leq \infty$, $0 < q \leq \infty$, $0 < s < \infty$. The space $B^s_{p,q}$ is defined as the set of all the functions $f \in \mathbb{L}_p$ such that:

$$\|f\|_{B^s_{p,q}} = \|f\|_p + \Big(\sum_{n=1}^{\infty}(n^s E_n(f,p))^q \frac{1}{n}\Big)^{1/q} < \infty, \quad \text{if } q < \infty,$$

$$(resp.\ \|f\|_{B^s_{p,\infty}} = \|f\|_p + \sup_{n\geq 1} n^s E_n(f,p) < \infty)$$

As $n \to E_n(f,p)$ is not increasing, we have an equivalent norm $\|f\|_{B^s_{p,q}}$ :

$$\|f\|_{B^s_{p,q}} \sim \|f\|_p + \Big(\sum_{j=0}^{\infty}(2^{js}E_{2^j}(f,p))^q\Big)^{1/q}, \quad \text{if } q < \infty,$$

$$(resp.\ \|f\|_{B^s_{p,\infty}} = \|f\|_p + \sup_{j\geq 0} 2^{js}E_{2^j}(f,p))$$

**Theorem 8.** *Let $f \in \mathbb{L}_p(\gamma_{\alpha,\beta})$ and*

$$\Lambda_j(f) = \sum_{\nu=1}^{2^j}\langle f, \psi_{j,\eta_\nu}\rangle \psi_{j,\eta_\nu}$$

*then*

1.

$$\|f\|_{B^s_{p,q}} \sim \|f\|_p + \Big(\sum_{j=0}^{\infty}(2^{js}\|\Lambda_j(f)\|_p)^q\Big)^{1/q} \tag{58}$$

   *with the usual modification for $q = \infty$.*

2.

$$\text{If } p < \infty, \quad \|f\|_{B^s_{p,q}} \sim \|f\|_p + \Big\{\sum_{j=0}^{\infty}(2^{js}(\sum_{\nu=1}^{2^j}|\langle f,\psi_{j,\eta_\nu}\rangle|^p\|\psi_{j,\eta_\nu}\|_p^p)^{1/p})^q\Big\}^{1/q} \tag{59}$$

   *with the usual modification for $q = \infty$.*

$$\|f\|_{B^s_{\infty,q}} \sim \|f\|_\infty + \Big\{\sum_{j=0}^{\infty}(2^{js}(\sup_\nu |\langle f,\psi_{j,\eta_\nu}\rangle|\|\psi_{j,\eta_\nu}\|_\infty))^q\Big\}^{1/q} \tag{60}$$

   *with the previous modification for $q = \infty$.*

*Proof.* If $Q \in \Pi_{2^{j-1}}$, $j \geq 1$, we have

$$\|\Lambda_j(f)\|_p = \|\Lambda_j(f-Q)\|_p \leq C_p\|f-Q\|_p$$

and hence

$$\|\Lambda_j(f)\|_p \leq E_{2^{j-1}}(f,p).$$

On the other hand
$$E_{2^j}(f,p) \le \sum_{m=j+1}^{\infty} \|\Lambda_m\|_p.$$

Therefore, if
$$\|\Lambda_m\|_p \le \varepsilon_m 2^{-ms}, \quad \varepsilon \in l_q(\mathbb{N}),$$

we have
$$E_{2^j}(f,p) \le \sum_{m=j+1}^{\infty} \varepsilon_m 2^{-ms} = 2^{-js} \sum_{m=j+1}^{\infty} \varepsilon_m 2^{-|m-j|s} = \delta_j 2^{-js},$$

with $\delta \in l_q(\mathbb{N})$, by a classical convolution result. Thus (58) is established.

Now, by (57) we have

$$\begin{aligned}
\|\Lambda_j(f)\|_p &\le C_p \Big(\sum_{\nu=1}^{2^j} |\langle f, \psi_{j,\eta_\nu}\rangle|^p \|\psi_{j,\eta_\nu}\|_{\mathbb{L}_p(\gamma_{\alpha,\beta})}^p\Big)^{1/p} \\
&\le C_p'(\|\Lambda_{j-1}f\|_p + \|\Lambda_j f\|_p + \|\Lambda_{j+1}f\|_p).
\end{aligned}$$

This combined with (58) readily implies (59) and (60). $\square$

**Remark 1.** *Suppose $\lambda_{j,\eta_\nu}$ is a family of numbers such that*

$$\Big\{\sum_{j=0}^{\infty} (2^{js} (\sum_{\nu=1}^{2^j} |\lambda_{j,\eta_\nu}|^p \|\psi_{j,\eta_\nu}\|_p^p)^{1/p})^q\Big\}^{1/q}$$

*(with the usual modifications for $p = \infty, q = \infty$)*
*and let*
$$f = c1 + \sum_{j=0}^{\infty} \sum_{\nu=1}^{2^j} \lambda_{j,\eta_\nu} \psi_{j,\eta_\nu},$$

*(which is, by Theorem 7, obviously defined in $\mathbb{L}_p$) $f \in B_{p,q}^s$ even though not necessarily*
$$\lambda_{j,\eta_\nu} = \langle f, \psi_{j,\eta_\nu}\rangle.$$

*This is due to Theorem 7, and since*
$$\langle f, \psi_{j,\eta_\nu}\rangle. = \langle A_{j-1} + A_j + A_{j+1}, \psi_{j,\eta_\nu}\rangle.$$

*where*
$$A_j = \sum_{\nu=1}^{2^j} \lambda_{j,\eta_\nu} \psi_{j,\eta_\nu}, \quad \|A_j\|_p \le C_p \Big(\sum_{\nu=1}^{2^j} \|\lambda_{j,\eta_\nu} \psi_{j,\eta_\nu}\|_p^p\Big)^{1/p}$$